%% file: main.tex
\documentclass[fleqn,11pt,letterpaper]{amsproc}
\usepackage{geometry}
\usepackage[T1]{fontenc}
\usepackage[utf8]{inputenc}
\usepackage{amssymb}
\usepackage{tikz-cd}
\include{math-shortcuts}
\renewcommand{\epsilon}{\varepsilon}
\usepackage{enumitem}
\usepackage{diagbox}
\tikzcdset{arrow style=tikz, diagrams={>=stealth}} 
\tikzcdset{every label/.append style = {font = \small}}
\DeclareMathOperator{\is}{is}
\let\SL\relax
\DeclareMathOperator{\SL}{\mathbf{SL}}
\let\SU\relax
\DeclareMathOperator{\SU}{\mathbf{SU}}
\let\GL\relax
\DeclareMathOperator{\GL}{\mathbf{GL}}
\let\SO\relax
\DeclareMathOperator{\SO}{\mathbf{SO}}
\let\Sp\relax
\DeclareMathOperator{\Sp}{\mathbf{Sp}}
\let\Lie\relax
\DeclareMathOperator{\Lie}{Lie}
\newcommand{\titlemath}[2]{\texorpdfstring{#1}{#2}}

\newcommand{\<}{\langle}
\renewcommand{\>}{\rangle}
\renewcommand{\tilde}{\widetilde}

\newenvironment{smallpmatrix}{\left(\begin{smallmatrix}}{\end{smallmatrix}\right)}
\renewcommand{\bar}{\overline}
\DeclareMathOperator{\Mahler}{\cM}
\renewcommand{\iff}{\quad\Leftrightarrow\quad}
\DeclareMathOperator{\Res}{Res}
\renewcommand{\emph}[1]{\textit{\color{blue}#1}}
\newcommand{\G}{\mathbf{G}}
\renewcommand{\C}{\mathbb{C}}
\newcommand{\symmetrized}[1]{#1+{#1}^{-1}}
\newcommand{\antisymmetrized}[1]{#1-{#1}^{-1}}
\newcommand{\psymmetrized}[1]{(\symmetrized{#1})}
\newcommand{\pantisymmetrized}[1]{(\antisymmetrized{#1})}
\newcommand{\myfrac}[2]{{#1}/{#2}}
\newcommand{\0}{\mbox{ }0\mbox{ }}

\makeatletter
\newif\ifnobrackets
\renewcommand\@cite[2]{\ifnobrackets\else[\fi{#1\if@tempswa , #2\fi}\ifnobrackets\else]\fi\nobracketsfalse}
\newcommand\nbcite{\nobracketstrue\cite}
\makeatother

\title{Arithmetic Groups and the Lehmer Conjecture}
\author{Lam Pham}
\address{L.P.:\ Department of Mathematics, Brandeis University, 415 South Street, 02453 Waltham, United States of America \& The Einstein Institute of Mathematics, Edmond J. Safra Campus, Givat Ram, The Hebrew University of Jerusalem Jerusalem, 91904, Israel}
\email{lampham@brandeis.edu}
\thanks{L.P.\ acknowledges the support of the \href{https://zuckerman-scholars.org/ourprograms/postdoc-program/}{Zuckerman STEM Leadership Program} (as Zuckerman Postdoctoral Scholar) and of the  \href{https://erc.europa.eu/}{ERC} (under HomDyn 833423).}
\author{Fran\c{c}ois Thilmany}
\address{F.T.:\ Department of mathematics, KU Leuven, Celestijnenlaan 200B, 3001 Leuven, Belgium}
\email{francois.thilmany@kuleuven.be}
\thanks{F.T.\ acknowledges the support of the \href{https://www.frs-fnrs.be/en/}{FNRS} (under CR FC 4057) and of the \href{https://www.kuleuven.be/}{KU Leuven} (under PDM 19/145).}
\date{\today}
\numberwithin{equation}{subsection}

\begin{document}

\begin{abstract}
We generalize a result of Sury \cite{Sury1992} and prove that uniform discreteness of cocompact lattices in higher rank semisimple Lie groups (first conjectured by Margulis \cite{Margulis1991}) is equivalent to a weak form of Lehmer's conjecture. We include a short survey of related results and conjectures.
\end{abstract}

\maketitle

\section{Introduction and main result}

\subsection{Margulis' arithmeticity theorem and conjecture}

Let \(\bfG\) be a connected semisimple \(\R\)-group with \(\rank_{\R}(\bfG)\geq 2\). Then, \(G=\bfG(\R)\) is a semisimple Lie group and Margulis' celebrated \emph{arithmeticity theorem} \cite{Margulis1975} states that every irreducible lattice \(\Gamma\subset G\) is \emph{arithmetic} (see \S \ref{section:Margulis-arithmeticity} for precise definitions and statements).

The starting point of this paper is the following consequence of arithmeticity for irreducible non-cocompact lattices \cite[IX, (4.21-A)]{Margulis1991}:

\begin{theorem*}[Margulis]
Assume that $\bfG$ has no $\R$-anisotropic factor. There is a neighbourhood \(U \subset G\) of the identity such that for any irreducible non-cocompact lattice \(\Gamma \subset G\), the intersection \(U \cap \Gamma\) consists of unipotent elements.
\end{theorem*}

Margulis \cite[IX, (4.21-B)]{Margulis1991} then conjectured that an analogous statement would hold for cocompact lattices. In fact, Margulis indicates that this conjecture would follow from a weaker form of Lehmer's conjecture, which we recall below (see namely Conjecture \ref{weak-Lehmer-conjecture}).

\begin{conjecture}[Margulis]\label{Margulis-conjecture}
Let \(\G\) be a connected semisimple \(\R\)-group. Suppose \(\rank_{\R}(\G)\geq 2\). Then there exists a neighborhood \(U\subset \G(\R)\) of the identity such that for any irreducible cocompact lattice \(\Gamma \subset \G(\R)\), the intersection \(U\cap\Gamma\) consists of elements of finite order.
\end{conjecture}

For the purpose of this paper, it will be useful to work with the following statement.
\begin{statement}[Margulis' conjecture for a family \(\scrT\) of semisimple \(\R\)-groups]
For each \(\bfG \in \scrT\), there exists a neighborhood \(U \subset \bfG(\R)\) of the identity such that for any irreducible cocompact \emph{arithmetic} lattice \(\Gamma\) in \(\bfG(\R)\), the intersection \(U \cap \Gamma\) consists of elements of finite order.
\end{statement}
\noindent Note that in view of the arithmeticity theorem, `Margulis' conjecture for higher rank groups' is in this sense simply Conjecture \ref{Margulis-conjecture}. 

\begin{remark*}
An even stronger statement than Conjecture \ref{Margulis-conjecture} holds for \(p\)-adic groups. Namely, if \(k\) is a non-archimedean local field of characteristic zero and \(\bfG\) is an algebraic \(k\)-group (of dimension \(>0\)), then there is an open neighborhood \(U\subset \bfG(k)\) of the identity such that each nontrivial element \(h \in U\) generates a non-discrete subgroup \cite[IX, (3.5)]{Margulis1991}. In particular, every lattice in \(\bfG(k)\) must intersect \(U\) \emph{trivially}.

In contrast, in the real case one cannot omit elements of finite order from the statement of Conjecture \ref{Margulis-conjecture}. 
Indeed, in general one can construct a sequence \((\Gamma_m)_{m \in \N}\) of cocompact arithmetic lattices in \(\bfG(\R)\) and torsion elements \(\gamma_m \in \Gamma_m\) tending to \(1\). In the case of \(\SL_2(\R)\) and \(\SL_2(\C)\), the reader can consult \cite[\S 12.5]{MaclachlanReid2003} for a thorough discussion on the existence of elements of arbitrary order in cocompact lattices. 
\end{remark*}

\subsection{Lehmer's conjecture} \label{Lehmer-subsec}

Let \(P\in\C[X]\) be a monic polynomial of degree \(d\) with roots \(\alpha_1,\,\hdots,\,\alpha_d \in \C\). The \emph{Mahler measure} of \(P\) is
\[
\Mahler(P)=\prod_{i=1}^d \max\{1,|\alpha_i|\}.
\]
In the following, let \(\alpha_1, \dots, \alpha_{s(P)}\) be an enumeration of the roots of \(P\) in \(\C\) which have absolute value strictly greater than 1, repeated according to their multiplicity, so that \(s(P)\) denotes their count (with multiplicity) and we may rewrite
\(
\Mahler(P) = \prod_{i=1}^{s(P)} |\alpha_i|.
\)
If \(\alpha \in \overline{\Q}\) is an algebraic integer, the Mahler measure \(\Mahler(\alpha)\) of \(\alpha\) will be defined as the Mahler measure \(\Mahler(P_\alpha)\) of its minimal polynomial \(P_\alpha\) over \(\Z\). The Mahler measure on algebraic integers is obviously invariant under the action of \(\Gal(\overline{\Q} / \Q)\).

The Mahler measure is multiplicative. By virtue of Kronecker's theorem, for \(P\) a monic, irreducible polynomial with integer coefficients, we have
\[
\Mahler(P)=1
\iff
P(X)=X,\ \text{or}\
P\
\text{is a cyclotomic polynomial.}
\]

In 1933, Lehmer \cite{Lehmer1933} asked whether one could find irreducible polynomials with integer coefficients whose Mahler measure gets arbitrarily close to 1 (but is not 1). It is conjectured that this is not possible:

\begin{conjecture}[Lehmer\protect{\footnote{Lehmer always insisted that he had not formulated his problem as a conjecture, although we will customarily refer to it as such.}}]\label{Lehmer-conjecture}
There exists \(\epsilon>0\) such that for any (irreducible) monic polynomial \(P\) with integer coefficients, either
\[
\Mahler(P)=1 \quad \text{or} \quad \Mahler(P)>1+\epsilon.
\]
\end{conjecture}

In fact, Lehmer's polynomial
\[
P_{\mathrm{Lehmer}} = X^{10}+X^9-X^7-X^6-X^5-X^4-X^3+X+1,
\]
for which \(\Mahler(P_{\mathrm{Lehmer}})=1.17628...\), is suspected to attain the smallest Mahler measure greater than 1.
Partial results towards Lehmer's conjecture are known; we list some of them below. \medskip

We will be concerned with the following weaker version of Lehmer's conjecture.

\begin{conjecture}[weak Lehmer]\label{weak-Lehmer-conjecture}
For each \(s \in \N\), there exists \(\epsilon(s)>0\) such that for any (irreducible) monic polynomial \(P\) with integer coefficients and \(s(P) \leq s\), either
\[
\Mahler(P) = 1 \quad \textrm{or} \quad \Mahler(P)>1+\epsilon(s).
\]
\end{conjecture}

For a given \(s\), we will call the statement in Conjecture \ref{weak-Lehmer-conjecture} \emph{Lehmer's conjecture at level \(s\)}. In this way, Conjecture \ref{weak-Lehmer-conjecture} could be described as ``Lehmer's conjecture at all levels'', and Lehmer's conjecture \ref{Lehmer-conjecture} as ``Lehmer's conjecture at all levels uniformly''.

\subsection{Main result}\label{main-theorem}

Fix an absolutely (almost) simple isotropic \(\R\)-group \(\bfF\) and consider, for each integer \(s \geq 1\), the family of semisimple \(\R\)-groups
\begin{equation*}
\scrT_\bfF^{(s)} = \left \{\prod_{i=1}^r \bfF \times \prod_{i=1}^t \Res_{\C/\R}(\bfF) \ \middle | \ r, t \in \N, \; r+2t \leq s \right \}. \label{eq:special-form}
\end{equation*}

The following theorem is the main result of this paper.

\begin{theorem*}
Let \(s \geq 1\). Then Margulis' conjecture for any of the families \(\scrT_\bfF^{(s)}\) defined above\footnote{As the statement in Margulis' conjecture is insensitive to isogenies (see namely the first paragraph of \S\ref{first-implication-subsec}), in the theorem, one could of course replace each element of the family \(\scrT_\bfF^{(s)}\) by another semisimple group isogenous to it.} implies Lehmer's conjecture at level \(s\).
In consequence, Margulis' conjecture  \ref{Margulis-conjecture} is equivalent to the weak version \ref{weak-Lehmer-conjecture} of Lehmer's conjecture.
\end{theorem*}

In view of the theorem, let us emphasize that Margulis' conjecture for products of simple groups of \emph{any fixed \(\R\)-type} and their extensions to \(\C\) is already sufficient to imply the weak Lehmer conjecture, which in turn implies Margulis' conjecture for \emph{all} semisimple groups (see \S \ref{section:Margulis-arithmeticity}).
For the sake of clarity, we will first carry the argument out in detail for type \(\sfA_n\) split in \S\ref{section:proof-of-main-theorem}. To complete the proof of the theorem, we then give a comprehensive treatment of all other \(\R\)-types in \S\ref{othergroups-sec}. It is noteworthy that the argument works for every \(\R\)-type in a similar way.

\medskip

It would be interesting to know whether Margulis' conjecture has implications for the full Lehmer conjecture \ref{Lehmer-conjecture}, beyond Conjecture \ref{weak-Lehmer-conjecture}. 
Connections between the full Lehmer conjecture and linear groups have already been brought to light, for example in the work of Breuillard \cite{Breuillard2007} and of Breuillard--Varju \cite{BreuillardVarju2020}. 

\subsection{Notation and conventions}
Throughout the paper, we adopt the conventions of Bourbaki. In particular, \(\N = \{0, 1, 2, \dots\}\) and we denote \(\N^* = \N \backslash \{0\}\).

\begin{itemize}
\item If \(P \in \C[x]\), \(s(P)\) denotes the number of roots of \(P\) in \(\C\) which have absolute value \(> 1\), counted with multiplicity.
\item \(\Mahler(P)\) denotes the Mahler measure of \(P \in \C[x]\).
\item If \(\alpha \in \overline{\Q}\) is an algebraic integer, \(P_\alpha\) denotes its minimal polynomial over \(\Z\).
\item \(\cO_K\) denotes the ring of integers of a number field \(K\).
\item \(\rM_n(A)\) denotes the set of \(n \times n\) matrices with entries in an algebraic structure \(A\), endowed with whichever structure is inherited from \(A\). 
\end{itemize}

\subsection*{Acknowledgements}

The authors would like to thank the anonymous referee for his helpful comments, which improved the clarity of the paper and led to a more thorough treatment of the exceptional groups.
The authors also thank E. Breuillard and G. Margulis for interesting conversations.

\section{A short history}

In this section, we provide a very short (and incomplete) survey of known results and references about the arithmetic and geometry of the famous Lehmer problem. For a more extensive treatment, we refer the reader to the surveys of Smyth \cite{Smyth2008,Smyth2015}, and Ghate and Hironaka \cite{GhateHironaka2001}.

\subsection{Some known results about the Lehmer conjecture}\label{Mahler-bounds}
Let \(\alpha \in \overline{\Q}\) be an algebraic integer which is not an integer nor a root of unity, and let \(d \geq 2\) denote the degree of \(P_\alpha\). Then, the following lower bound for \(M(\alpha)\) is known:
\begin{equation}\label{eq:Voutier}
\Mahler(\alpha)> 1+\frac{1}{4}\left(\frac{\log\log d}{\log d}\right)^3.
\end{equation}
This was proved by Voutier \cite[Theorem]{Voutier1996}, improving the bound of Dobrowolski \cite[Theorem 1]{Dobrowolski1979} (who obtained a factor of \(1/1200\) instead of \(1/4\)). Another bound due to Laurent\footnote{In a private communication, Michel Laurent informed us that this bound was not published and that Schinzel \cite{Schinzel1973} proved the sharp inequality \(\Mahler(\alpha)\geq \big(\tfrac{1+\sqrt{5}}{2}\big)^{d/2}\), for every \(\alpha\) which is a totally real algebraic integer of degree \(d\geq 2\).} (1983), emphasizing the number of real roots, is given by Margulis \cite[p.~322]{Margulis1991}: if \(P\) is a non-cyclotomic polynomial of degree \(d\geq 2\) with \(r\) real roots, then
\[
\Mahler(\alpha)\geq c^{r^2/d \log(1+\frac{d}{r})},
\]
where \(c > 0\) is an absolute constant.

\medskip

Recall that a polynomial \(P \in \Z[X]\) is called \emph{palindromic} if
\[
P(X)=X^{\deg(P)} P(X^{-1}).
\]
Every palindromic polynomial \(P\) of odd degree is divisible by \(X+1\), hence irreducible palindromic polynomials of degree \(>1\) must have even degree.

Smyth \cite{Smyth1971} proved that the polynomial \(P_{\mathrm{Smyth}}(X)=X^3-X-1\) had the smallest Mahler measure among non-palindromic polynomials.

\begin{theorem*}[\cite{Smyth1971}]
Let \(\alpha\in\bar{\Q}\). If \(\Mahler(P_\alpha) < \Mahler(P_{\mathrm{Smyth}})\), then \(P_\alpha\) is palindromic.
\end{theorem*}
\noindent This effectively reduces Conjectures \ref{Lehmer-conjecture} and \ref{weak-Lehmer-conjecture} to palindromic polynomials.

Although we will not use it, we note in passing another instance for which the Lehmer conjecture is known to hold. 
Borwein, Dobrowolski, and Mossinghoff \cite{BorweinDobrowolskiMossinghoff2007} proved that if \(P\) is a polynomial of degree \(d\) without cyclotomic factors, all of whose coefficients are odd integers, then
\[
\log \Mahler(P)\geq \frac{\log 5}{4} \left(1-\frac{1}{d+1}\right).
\]

\subsection{The case of Salem numbers}

Recall that a \emph{Salem number} is an algebraic integer \(u\in \R\) which is \(>1\) and all of whose Galois conjugates in \(\C\) have absolute value \(\leq 1\), with at least one of absolute value \(=1\).
Hence, any Salem number \(\alpha\) satisfies \(s(P_\alpha)=1\) where \(P_\alpha\) is the minimal polynomial of \(\alpha\), and it is clear that \(P_\alpha\) is palindromic. Conversely, any irreducible palindromic polynomial \(P \in \Z[x]\) with \(s(P)=1\) and \(\deg P \geq 4\) is the minimal polynomial of a Salem number. The following conjecture is thus equivalent to Lehmer's conjecture at level \(s=1\) (see the paragraph after Conjecture \ref{weak-Lehmer-conjecture}).

\begin{conjecture}[Salem, arithmetic version]\label{Salem-conjecture}
There exists \(\epsilon>0\) such that every Salem number \(\alpha\) satisfies \(\alpha >1+\epsilon\).
\end{conjecture}

Salem numbers constitute an important family of algebraic numbers. For example, the polynomial \(P_{\text{Lehmer}}\) with the smallest known Mahler measure (see \S \ref{Lehmer-subsec}) turns out to be the minimal polynomial of a {Salem number}.

\medskip

A beautiful connection with discrete subgroups of Lie groups was established by Sury \cite{Sury1992}: he proved that Conjecture \ref{Salem-conjecture} was equivalent to the following conjecture.

\begin{conjecture}[Salem, geometric version]\label{Sury-conjecture}
There exists a neighborhood \(U\subset \SL_2(\R)\) of the identity such that for any torsion-free cocompact arithmetic lattice \(\Gamma\subset \SL_2(\R)\), we have \(\Gamma\cap U=\{e\}\).
\end{conjecture}

\begin{remark*}
Conjecture \ref{Sury-conjecture} is simply Margulis' conjecture for \(\SL_2\), hence Theorem \ref{main-theorem} for \(s=1\) includes Sury's result as a special case. 
In fact, Theorem \ref{main-theorem} for \(s=1\) shows that Margulis' conjecture \ref{Margulis-conjecture} for any one of the isotropic absolutely (almost) simple \(\R\)-groups (e.g.\ \(\SL_{n}\), \(\SO_q\), etc.) implies Conjecture \ref{Salem-conjecture}.
\end{remark*}

\subsection{Lengths of shortest geodesics in arithmetic hyperbolic orbifolds}

Extending the picture to Kleinian groups, Neumann and Reid \cite{NeumannReid1992} formulated the following conjecture.

\begin{conjecture}[Short geodesic conjecture]\label{short-geodesic-conjecture}
There is a positive universal lower bound for the lengths of closed geodesics in arithmetic hyperbolic \(2\)- and \(3\)-orbifolds.
\end{conjecture}

The short geodesic conjecture for hyperbolic 2-orbifolds is equivalent to Conjecture \ref{Sury-conjecture}, hence to \ref{Salem-conjecture}. On the other hand, the short geodesic conjecture for hyperbolic 3-orbifolds is equivalent to the following \emph{complex Salem conjecture}. Call an algebraic integer \(\alpha \in \C\) a \emph{complex Salem number} if \(\alpha\) is not real, \(\alpha\) and its complex conjugate \(\overline{\alpha}\) are the only two Galois-conjugates of \(\alpha\) in \(\C\) of absolute value \(> 1\), and \(\alpha\) has at least one conjugate of absolute value \(=1\).

\begin{conjecture}[Complex Salem] \label{complex-Salem-conjecture}
There exists \(\epsilon > 0\) such that every complex Salem number \(\alpha\) satisfies \(|\alpha| > 1 + \epsilon\).
\end{conjecture}

In fact, the short geodesic conjecture for hyperbolic \(3\)-orbifolds implies the conjecture for \(2\)-orbifolds. This can be seen arithmetically, as the complex Salem conjecture implies the traditional Salem conjecture. Indeed, if \(\alpha\) is a Salem number, then  \(P_\alpha(-x^2)\) is the minimal polynomial of a complex Salem number with the same Mahler measure as \(\alpha\).

For a detailed introduction to Conjecture \ref{short-geodesic-conjecture}, we refer the reader to the book of Maclachlan and Reid \cite{MaclachlanReid2003}.

\begin{remark*}
Similarly as for Conjecture \ref{Salem-conjecture}, the proof of Theorem \ref{main-theorem} (for \(s=2\), \(r=0\)) shows that Margulis' conjecture for the restriction to \(\R\) of any simple \(\C\)-group implies Conjecture \ref{complex-Salem-conjecture}. 
\end{remark*}

\medskip

Salem numbers can also be used to obtain a lower bound for the length of closed geodesics in {noncompact arithmetic hyperbolic orbifolds of \emph{even} dimension \(n\)}. This was done by Emery, Ratcliffe and Tschantz \cite{EmeryRatcliffeTschantz2019}. More precisely, for any integer \(n\geq 2\), let \(\cH^n\) denote hyperbolic \(n\)-space, and define
\[
\beta_n = \min\{\log \alpha \mid \alpha~\text{is a Salem number with}~\deg P_\alpha \leq  n\}.
\]
For any even dimension \(n\), if \(\Gamma\subset\Isom(\cH^n)\) is a non-uniform arithmetic lattice, then the length of any closed geodesic in \(\cH^n/\Gamma\) is at least \(\beta_n\) \cite[Corollary 1.3]{EmeryRatcliffeTschantz2019}. It follows that Conjecture \ref{Salem-conjecture} is equivalent to the existence of a uniform lower bound on the length of closed geodesics in non-compact arithmetic hyperbolic orbifolds of even dimension.

\subsection{Homotopy type of locally symmetric spaces}

Let us conclude this section with a consequence of Conjecture \ref{Margulis-conjecture} in the context of locally symmetric spaces. Let \(S\) be a symmetric space (e.g., \(S=G/K\) where \(G\) is a semisimple Lie group, \(K\subset G\) a maximal compact subgroup). Following Gelander \cite{Gelander2004}, an \emph{\(S\)-manifold} is a complete Riemannian manifold locally isometric to \(S\), i.e., a manifold of the form \(\Gamma\backslash S\), where \(\Gamma\subset \Isom(S)\) is a discrete torsion-free subgroup. It is irreducible if \(\Gamma\) is an irreducible lattice. For \(d,v\in\N\), a \emph{\((d,v)\)-simplicial complex} is a simplicial complex with at most \(v\) vertices, all of degree at most \(d\). An interesting consequence of Conjecture \ref{Margulis-conjecture} is the following conjecture of Gelander.

\begin{conjecture}[\cite{Gelander2004}]\label{Gelander-conjecture}
For any symmetric space \(S\) of noncompact type, there are constants \(\alpha(S)\), \(d(S)\), such that any irreducible \(S\)-manifold \(M\) (assumed to be  arithmetic if \(\dim(S) = 3\)) is homotopically equivalent to a \((d(S), \alpha(S) \vol(M))\)-simplicial complex.
\end{conjecture}

Conjecture \ref{Gelander-conjecture} was recently proved by Fraczyk \cite[Theorem 1.16]{Fraczyk2021} for arithmetic 3-manifolds.

\section{Margulis' arithmeticity theorem} \label{section:Margulis-arithmeticity}

In this section, we briefly review Margulis' arithmeticity results \cite[Chapter IX]{Margulis1991} in a form that will be used to show that Conjecture \ref{weak-Lehmer-conjecture} implies Conjecture \ref{Margulis-conjecture}. Let us first recall the definitions.

Let \(\bfG\) be a connected semisimple \(\R\)-group and let \(\prod_{i\in I} \bfG_i\) be a decomposition of \(\bfG\) as an almost direct product of almost \(\R\)-simple \(\R\)-subgroups. 
For any subset \(J\subset I\), we write \(\bfG_J=\prod_{i\in J} \bfG_i\). 
We denote by \(\bfG^{\is}\) (resp.\ \(\bfG^\mathrm{anis}\)) the subgroup of \(\bfG\) which is the almost direct product of the \(\R\)-isotropic (resp.\ \(\R\)-anisotropic) factors of \(\bfG\). We also let \(G=\bfG(\R)\), \(G_J=\bfG_J(\R)\) and \(G^{\is}=\bfG^{\is}(\R)\).

A lattice \(\Gamma\subset G\) is said to be \emph{irreducible} if for any non-empty proper subset \(J\subset I\), the index \([\Gamma:(\Gamma\cap G_J) \cdot (\Gamma\cap G_{I\setminus J})]\) is infinite.

An irreducible lattice \(\Gamma\subset G\) such that \(G^{\is}\cdot\Gamma\) is dense in \(G\) is called \emph{arithmetic} if there exist a connected non-commutative almost \(\Q\)-simple \(\Q\)-group \(\bfH\) (endowed with some \(\Z\)-structure), and an \(\R\)-epimorphism \(\tau:\bfH\to\bfG\) such that:
\begin{enumerate}[label=(\roman*)]
\item the Lie group \((\ker\tau)(\R)\) is compact;
\item the subgroups \(\tau(\bfH(\Z))\) and \(\Gamma\) are commensurable.
\end{enumerate}

\subsection{The arithmeticity theorem} \label{subsection:Margulis-arithmeticity}

We can now state Margulis' celebrated arithmeticity theorem.

\begin{theorem*}[{\cite[IX, (1.16)]{Margulis1991}}]
Let \(\bfG\) be a connected semisimple \(\R\)-group and \(\Gamma\) an irreducible lattice in \(\bfG(\R)\), with \(\bfG^{\is}(\R)\cdot\Gamma\) dense in \(\bfG(\R)\). Suppose that \(\rank_{\R} \bfG \geq 2\). Then the lattice \(\Gamma\) is arithmetic.
\end{theorem*}

In the remainder of this section, we will assume that \(\rank_{\R}(\bfG)\geq 2\), that \(G\) has no compact factors (i.e.\ \(\bfG^{\is} =  \bfG\)) and has trivial center. Let \(\Gamma \subset G\) be an irreducible lattice. Then the following facts are among the key steps of the proof of the arithmeticity theorem. We refer the reader to \cite[\S 6.1]{Zimmer1984} or \cite[\S11.5]{Benoist2008} for proofs and details.

\begin{enumerate}[leftmargin=\mathindent,%
label=(\thesubsection.\arabic*),%
labelwidth=\mathindent,%
labelsep=0pt,
parsep=4pt,%
align=left,%
ref=\thesubsection.\arabic*]
\item The trace field \(K=\Q(\tr\Ad\Gamma)\) of \(\Gamma\), the field generated by the set \(\{\tr(\Ad(\gamma))\,|\, \gamma\in\Gamma\}\), is a number field. This follows from Margulis' superrigidity theorem, together with the fact that \(\Gamma\) is finitely generated.

\item Since \(\Gamma\) is Zariski-dense in \(\bfG(\R)\), \(\bfG\) can be defined over \(K\). That is, there is a \(K\)-group \(\bfG'\) and a place \(v_0\) of \(K\) such that \(K_{v_0} = \R\), \(\bfG \cong \bfG'\) as \(\R\)-groups, and the image of \(\Gamma\) under this isomorphism lies in \(\bfG'(K)\). In the following, we will identify \(\bfG\) and \(\bfG'\), and simply assume that \(\Gamma \subset \bfG(K)\).

\item
\label{fancyarithmeticity-point}
There is a semisimple \(\Q\)-group \(\bfH\) (endowed with a \(\Z\)-structure) with trivial center, an \(\R\)-epimorphism \(\tau: \bfH \to \bfG\) with \((\ker \tau) (\R)\) compact, and a homomorphism \(\iota: \Gamma \to \bfH(\Q)\) such that \(\tau \circ \iota\) is the identity, \(\iota(\Gamma)\) is Zariski dense in \(\bfH\) and \(\iota(\Gamma)\) is commensurable with \(\bfH(\Z)\).
The group \(\bfH\) can be constructed as the restriction of scalars \(\Res_{K/\Q}(\bfG)\) of \(\bfG\) from \(K\) to \(\Q\). 
\end{enumerate}

\subsection{The weak Lehmer conjecture implies Margulis' conjecture} \label{first-implication-subsec}
We start by indicating how Conjecture \ref{weak-Lehmer-conjecture} implies Conjecture \ref{Margulis-conjecture}. Note that as long as the lattice \(\Gamma\) is arithmetic, the argument given below works equally well for groups of rank one. It shows in fact that Conjecture \ref{weak-Lehmer-conjecture} implies Margulis' conjecture for all semisimple \(\R\)-groups (in the sense defined below the statement of Conjecture \ref{Margulis-conjecture}). 
Let thus \(\bfG\) be a semisimple group and \(\Gamma\) an irreducible arithmetic lattice in \(\bfG(\R)\). 

Without loss of generality, we may assume that \(\bfG\) has trivial center and is without anisotropic factors.
Indeed, if \(\bfG\) has center \(\bfC\), and \(U'\) is a neighborhood of \((\bfG/\bfC)(\R)\) as in \ref{Margulis-conjecture}, then the preimage \(U\) of \(U'\) under the canonical map \(\pi: \bfG(\R) \to (\bfG/\bfC)(\R)\) has the required property: if \(\Gamma\) is an irreducible cocompact lattice in \(\bfG(\R)\), then \(\pi(\Gamma)\) is an irreducible cocompact lattice in \((\bfG/\bfC)(\R)\); hence if \(\gamma \in U\), \(\pi(\gamma)\) must have finite order, and since \(\bfC(\R)\) is finite, so does \(\gamma\).
Similarly, if \(\bfG = \bfG^{\is} \times \bfG^{\mathrm{anis}}\) and \(U'\) is a neighborhood of \(\bfG^{\is}\) as in \ref{Margulis-conjecture}, then the preimage \(U\) of \(U'\) under the canonical map \(\pi: \bfG(\R) \to \bfG^{\is}(\R)\) has the required property: if \(\Gamma\) is an irreducible cocompact lattice in \(\bfG(\R)\), then \(\pi(\Gamma)\) is an irreducible cocompact lattice in \(\bfG^{\is}(\R)\); hence if \(\gamma \in U\), \(\pi(\gamma)\) has finite order \(m\), and \(\gamma^m \in \bfG^{\mathrm{anis}}(\R) \cap \Gamma\). As \(\bfG^{\mathrm{anis}}(\R)\) is compact and \(\Gamma\) is discrete, the latter is a finite group and \(\gamma\) has finite order.

Let \(\bfH\) be the group obtained in \eqref{fancyarithmeticity-point}, so that we have the following diagram.
\[
\begin{tikzcd}[row sep=1.2cm, column sep=1.2cm]
& \Gamma \arrow[d, hook] \arrow[ld, "\iota"'] \\
\bfH(\R) \arrow[r, "\tau"',swap] \arrow[d, "\Ad"'] & \bfG(\R) \arrow[d, "\Ad"]\\
\mathrm{GL}(\Lie(\bfH(\R))) \arrow[r] &  \mathrm{GL}(\Lie(\bfG(\R)))
\end{tikzcd}
\]
Recall also from \eqref{fancyarithmeticity-point} that \(\tau\) has compact kernel and that \(\iota(\Gamma)\) is commensurable with \(\bfH(\Z)\).

Since the adjoint representation \(\bfH \to \Ad \bfH\) is defined over \(\Q\), we can find a finite-index subgroup \(\Lambda\) of \(\bfH(\Z)\) for which \(\Ad(\Lambda) \subset \Ad(\bfH)(\Z)\) (see for example \cite[I, (3.1.1)]{Margulis1991}); in particular, \(\Ad(\Lambda)\) preserves a lattice in \(\Lie(\bfH(\R))\). Since \(\Lambda\) and \(\iota(\Gamma)\) are commensurable, \(\Ad(\iota(\Gamma))\) also stabilizes a lattice in \(\Lie(\bfH(\R))\) \cite[IX, (4.19)]{Margulis1991}. Hence the characteristic polynomials of the elements of \(\Ad(\iota(\Gamma))\) have integer coefficients.

Let us write \(\bfH(\R) = F \times K\) as a direct product where \(K\) is compact and \(F\) is without compact factors. The morphism \(\tau\) then induces an isogeny \(F \to \bfG(\R)\) and \(\ud \tau\) restricts to an isomorphism \(\Lie(F) \to \Lie(\bfG(\R))\). Let \(x \in \bfH(\R)\) and write \(\Ad_{\bfH(\R)}(x)=\Ad_F(y) \oplus \Ad_K(z)\) for some \(y \in F\), \(z \in K\). If we denote \(P_x\), \(P_y\), \(P_z\) the respective characteristic polynomials of \(\Ad_{\bfH(\R)}(x)\), \(\Ad_F(y)\), \(\Ad_K(z)\), we have that \(P_x = P_y \cdot P_z\).
Because \(K\) is compact, all roots of \(P_z\) in \(\C\) must have absolute value 1; in consequence, \(\Mahler(P_z) = 1\), hence \(\Mahler(P_x) = \Mahler(P_y)\).
Moreover, \(P_x\) satisfies \(s(P_x) = s(P_y) \leq \dim F = \dim \bfG\).

Now pick \(\gamma \in \Gamma\) and apply the last paragraph to \(x = \iota(\gamma)\). Since \(\tau(\iota(\gamma))=\gamma\) and \(\ud \tau \circ \Ad(\iota(\gamma)) = \Ad(\gamma) \circ \ud \tau\),
we obtain that the characteristic polynomial \(P_{\gamma}\) of \(\Ad_{\bfG(\R)}(\gamma)\) equals \(P_{y}\). We have in turn
\[
\Mahler(P_\gamma) = \Mahler(P_{\iota(\gamma)})
\quad\text{and}\quad
s(P_\gamma) = s(P_{\iota(\gamma)}).
\]

Let \(f:\bfG(\R) \to [1, \infty[\) be defined by \(f(g) = \Mahler(P_g)\), where as before \(P_g\) denotes the characteristic polynomial of \(\Ad_{\bfG(\R)}(g)\); note that \(f\) is a continuous function.
If the weak version of Lehmer's conjecture holds at level \(s= \dim \bfG\), we can find \(\epsilon > 0\) such that any polynomial \(P\) with integer coefficients and \(s(P) \leq \dim \bfG\) satisfies either \(\Mahler(P) = 1\) or \(\Mahler(P) > 1+\epsilon\).
This applies to \(P_{\iota(\gamma)}\) for any \(\gamma \in \Gamma\): as we observed, \(P_{\iota(\gamma)}\) has integer coefficients and \(s(P_{\iota(\gamma)}) \leq \dim \bfG\).
Thus, by the above, the open neighborhood \(U = f^{-1}([1, \epsilon[)\) of \(1\) in \(\bfG(\R)\) is such that \(U \cap \Gamma\) consists of elements \(\gamma\) for which \(\Mahler(P_\gamma) = \Mahler(P_{\iota(\gamma)}) = 1\). This means that for \(\gamma \in U\), \(P_\gamma\) is a product of cyclotomic polynomials. Since \(\Gamma\) is cocompact, \(\Ad(\gamma)\) is semisimple and thus \(\Ad(\gamma)\), hence also \(\gamma\), have finite order.

\section[Proof of the main theorem for \(\SL_n\)]{Proof of the main theorem for \(\mathrm{SL}_n\)}
\label{section:proof-of-main-theorem}

We now proceed to prove Theorem \ref{main-theorem} (in particular, that Conjecture \ref{Margulis-conjecture} implies Conjecture \ref{weak-Lehmer-conjecture}). To this end, we will assume Conjecture \ref{weak-Lehmer-conjecture} fails and construct a sequence of cocompact lattices in suitable groups of the family \(\scrT_\bfF^{(s)}\) violating Margulis' conjecture for this family.
As it is more transparent, we first give the full argument for \(\bfF = \SL_{n}\), i.e.\ for the family
\begin{equation*}\label{eq:special-form-SL}
\scrT_{\SL_{n}}^{(s)} = \left \{\prod_{i=1}^r \SL_{n} \times \prod_{i=1}^t \Res_{\C/\R}(\SL_{n}) \ \middle | \ r, t \in \N, \; r+2t \leq s \right \}.
\end{equation*}
This is already sufficient to establish the equivalence between Conjectures \ref{Margulis-conjecture} and \ref{weak-Lehmer-conjecture}. 
In the last section (\S\ref{othergroups-sec}), we then indicate the corresponding changes for \(\bfF\) an arbitrary absolutely (almost) simple isotropic \(\R\)-group, completing the proof of Theorem \ref{main-theorem}. 

\subsection{Reduction to palindromic polynomials with control on the archimedean places} \label{reduction-subsec}
Given \(P \in \Z[x]\), let us denote as before \(\alpha_1, \dots \alpha_{s(P)}\) the roots of \(P\) in \(\C\) of absolute value \(>1\), and label them in such a way \(\alpha_1, \dots, \alpha_{r(P)}\) lie in \(\R\) and \(\alpha_{r(P)+1}, \dots, \alpha_{s(P)}\) do not. In addition, let us order the latter roots so that \(\overline{\alpha}_{r(P)+i} = \alpha_{(r(P)+s(P))/2 + i}\) for \(1 \leq i \leq (s(P)-r(P))/2\). For each \(s' \in \N\), we will consider the set of polynomials
\[
\scrP_{\leq s'}=
\big\{P\in\Z[X]\,|\, P\ \text{is monic, irreducible, palindromic, and } s(P) \leq s' \big\},
\]
and for \(r \leq s \leq s'\), its subsets
\[
\scrP_{s,r} =\big\{P \in \scrP_{\leq s'} \mid s(P) = s \textrm{ and } r(P) = r\big\}.
\]
By construction, \(\scrP_{\leq s'}\) is the disjoint union \(\bigcup_{s=0}^{s'} \bigcup_{r=0}^{s} \scrP_{s,r}\).

Suppose that Conjecture \ref{weak-Lehmer-conjecture} does not hold. That is, there are \(s' \in \N\) and a sequence \((P_m)_{m \in \N}\) of monic, irreducible polynomials with integer coefficients such that \(s(P_m) \leq s'\) and \(\Mahler(P_m) \to 1\) while \(\Mahler(P_m) > 1\). By virtue of Smyth's theorem (\S \ref{Mahler-bounds}), we may assume that each \(P_m\) is palindromic, i.e.\ that \(P_m \in \scrP_{\leq s'}\). Moreover, up to extracting an appropriate subsequence, we may assume that \(P_m \in \scrP_{s,r}\) for some fixed integers \(r \leq s \leq s'\).\footnote{If \(s'\) were to be the smallest integer for which Lehmer's conjecture at level \(s'\) fails, then obviously \(s=s'\).} Of course, \(s > 0\).
In view of the bound \eqref{eq:Voutier}, it must be that the sequence \((\deg P_m)_m\) is unbounded. We may thus additionally assume that \(\deg P_m > 2s\). This implies the following important feature: \(P_m\) must have a root of absolute value \(1\). Indeed, \(P_m\) has \(s\) roots of absolute value \(> 1\), and because \(P_m\) is palindromic, as many roots of absolute value \(<1\).

The discussion above shows that the negation of Lehmer's conjecture at level \(s'\) (cf.\ the paragraph following Conjecture \ref{weak-Lehmer-conjecture}) amounts to the following statement.

\begin{statement}[Negation of Lehmer's conjecture at level \(s'\)]
There are fixed integers \(r\in\N,\, s\in\N^*\) with \(r \leq s \leq s'\)
such that:
\begin{enumerate}[leftmargin=\mathindent,%
label=(L\(_{s,r}\)),%
labelwidth=\mathindent,%
labelsep=0pt,
parsep=4pt,%
align=left,%
ref={L\(_{s,r}\)}]
\item\label{L}
for any \(\epsilon > 0\), there is a polynomial \(P \in \scrP_{s,r}\) with at least one root in \(\C\) of absolute value \(1\), for which \(1< \Mahler(P) < 1 + \epsilon\).
\end{enumerate}
\end{statement}

Using statement \eqref{L} as the main ingredient, for each such pair \((s,r)\), we will construct in the semisimple \(\R\)-group
\begin{equation}\tag{\(\star\)}\label{eq:Gsr}
\bfG_{s,r} = \prod_{i=1}^{r} \SL_n \times \prod_{i=1}^{(s-r)/2} \Res_{\C/\R} (\SL_n) \qquad (n\geq 2),
\end{equation}
a sequence of cocompact lattices \(\Gamma_m < \bfG_{s,r}(\R)\) and a sequence of elements \(\gamma_m \in \Gamma_m\) of infinite order, such that \(\gamma_m \to e \in \bfG_{s,r}(\R)\) as \(m \to \infty\). This shall provide a counterexample to Margulis' conjecture for \(\bfG_{s,r}\), thus proving the first part of theorem \ref{main-theorem}.
The equivalence of Margulis' conjecture \ref{Margulis-conjecture} and the weak version \ref{weak-Lehmer-conjecture} of Lehmer's conjecture then follows immediately by combining this with \S\ref{first-implication-subsec}.

\subsection{The number fields \titlemath{\(K\)}{K} and \titlemath{\(L\)}{L}}
\label{numberfields-subsec}

Let \(r\leq s\) be as above. Pick a polynomial \(P \in \scrP_{s,r}\), set \(2d = \deg P\), and let \(L = \Q(\alpha)\) denote the number field generated over \(\Q\) by a root \(\alpha\) of \(P\). Since \(P\) is palindromic, \(\alpha^{-1}\) is also a root of \(P\). In particular, \(\alpha^{-1}\) is integral over \(\Z\), and the assignment \(\tau:\alpha \mapsto \alpha^{-1}\) defines a non-trivial automorphism of \(L\) (which restricts to an automorphism of its ring of integers \(\cO_L\)). Let \(K = \Q\psymmetrized{\alpha}\) denote the subfield of \(L\) generated over \(\Q\) by \(\symmetrized{\alpha}\). Since \(K\) is the fixed field of \(\tau\), \(L\) is a quadratic extension of \(K\) whose Galois group is \(\{\id, \tau\}\). Note that the minimal polynomial of \(\alpha\) over \(K\) is \(X^2-\psymmetrized{\alpha}X +1\).

For \(1 \leq i \leq s\), let \(\sigma_i\) denote the embedding \(L \to \C\) defined by \(\sigma_i(\alpha) = \alpha_i\). By definition, \(L\) has \(r\) real embeddings \(\sigma_1,\, \dots,\, \sigma_{r}\) and \(s - r\) complex embeddings \(\sigma_{r+1},\, \dots,\, \sigma_s\) for which the image of \(\alpha\) has absolute value \(>1\).
\(L\) also has \({d} - s\) pairs of conjugate complex embeddings \(\sigma_{s+1},\, \overline{\sigma}_{s+1},\, \dots,\, \sigma_{d},\, \overline{\sigma}_{d}\) for which \(\alpha\) maps to an element of absolute value \(1\). Observe that precomposition with \(\tau\) sends \(\{\sigma_1, \dots, \sigma_{r}\}\) and \(\{\sigma_{r+1}, \dots, \sigma_{s}\}\) respectively to the set of real and complex embeddings of \(L\) for which the image of \(\alpha\) has absolute value \(<1\). Similarly, \(\sigma_{i} \circ \tau = \overline{\sigma}_{i}\) for \(s+1 \leq i \leq d\).
By construction, the embeddings \(\sigma_i\) and \(\sigma_i \circ \tau\) for \(1 \leq i \leq d\) agree on \(K\). This shows that \(\{\sigma_{1}, \dots, \sigma_{d}\}\) is the complete set of embeddings of \(K\) into \(\C\) (we omit the restriction to \(K\) from the notation). The image \(\sigma_i\psymmetrized{\alpha}\) of the generator of \(K\) belongs to \(\R\) if and only if \(\sigma_i(\alpha) \in \R\) or \(|\sigma_i(\alpha)| = 1\). In consequence, restricted to \(K\), the embeddings \(\sigma_1, \dots, \sigma_r\) and \(\sigma_{s+1}, \dots, \sigma_{d}\) are real, \(\sigma_{r+1}, \dots, \sigma_{s}\) are complex (coming in conjugate pairs), and \(K\) has signature \((r - s + d, (s - r)/2)\).

\subsection{The \(K\)-group \(\bfG\)} \label{group-subsec}
Let \(h: L^n \times L^n \to L\) be the \(\tau\)-hermitian form given by
\[
h(x,y) = x_1 \tau(y_1) + \dots + x_n \tau(y_n) \qquad x,\,y\in L^n,
\]
and let \(\bfG = \SU_h\) be the special unitary group associated to the form \(h\). \(\bfG\) is a linear algebraic \(K\)-group, whose group of \(K\)-points is (isomorphic to) the group of \(
h\)-unitary matrices in \(\rM_n(L)\) of determinant \(1\).\footnote{In fact, since the equations defining \(\bfG\) can be taken with coefficients in \(\cO_K\), \(\bfG\) can also be viewed as an \(\cO_K\)-group scheme.}
The \(K \otimes_\Q \R\)-points of \(\bfG\) can be computed easily by studying the behavior of the extension \(L/K\) and of the form \(h\) under the different embeddings \(\{\sigma_1, \dots, \sigma_d\}\) of \(K\).

In the present setting, the extension \(L/K\) splits at the places \(\sigma_1, \dots, \sigma_r\) (since these extend to real places of \(L\)) and also at the places \(\sigma_{r+1}, \dots, \sigma_{s}\) (since these are complex places of \(K\)), but not at \(\sigma_{s+1}, \dots, \sigma_{d}\) (since there we have \(L \otimes_K K_{\sigma_i} \cong \C\)).
At these last places, the hermitian form \(h\) becomes the standard hermitian form on \(\C^n \times \C^n\) after identifying the completion of \(L/K\) with \(\C / \R\) via \(\sigma_i\) (\(s+1 \leq i \leq d\)).
Altogether, we have
\[
\bfG(K \otimes_\Q \R) \cong \prod_{i=1}^{r} \SL_n(\R) \times \prod_{i=1}^{(s-r)/2} \SL_n(\C) \times \prod_{i=1}^{d-s} \SU_n(\R),
\]
where \(\SU_n(\R)\) denotes the standard anisotropic special unitary group over \(\R\).

\subsection{The lattice \(\Gamma\) and the element \(\gamma\)}
\label{lattice-subsec}

Let \(\Gamma\) be the group of \(h\)-unitary \(n \times n\) matrices of determinant \(1\) with entries in \(\cO_L\).\footnote{With this definition, \(\Gamma\) is (commensurable to) the \(\cO_K\)-points of \(\bfG\) when it is viewed as an \(\cO_K\)-group scheme as above.}
A classical theorem of Borel and Harish--Chandra \cite{BorelHarish-Chandra1962} states that \(\Gamma\) is a lattice in \(\bfG(K \otimes_\Q \R)\) when embedded diagonally using the inequivalent archimedean places of \(K\) among \(\sigma_1, \dots, \sigma_d\). (Apply \cite[Theorem 7.8]{BorelHarish-Chandra1962} to \(\Res_{K/\Q}(\bfG)\) and observe that \(\Res_{K/\Q}(\bfG)(\R) \cong \bfG(K \otimes_\Q \R)\); under this isomorphism, \(\Res_{K/\Q}(\bfG)(\Z)\) is commensurable to \(\Gamma\).) By virtue of the strong approximation theorem (for number fields), \(\Gamma\) is an irreducible lattice.
Moreover, if \(P\) has at least one root in \(\C\) of absolute value \(1\) (i.e.\ if \(s < d\)), we claim that \(\Gamma\) is a cocompact lattice in \(\bfG(K \otimes_\Q \R)\). Indeed, for any embedding \(\sigma: K \to \C\) which extends to \(L\) in such a way \(|\sigma(\alpha)|=1\), the group \(\bfG(K_\sigma) \cong \SU_n(\R)\) is compact, and hence \(\bfG\) must be \(K\)-anisotropic (as it is anisotropic over \(K_\sigma\)). The claim then follows from Godement's criterion \cite[Theorem 11.8]{BorelHarish-Chandra1962}.

Let \(\bfG_{s,r}\) be the \(\R\)-group defined in~ \ref{reduction-subsec}\eqref{eq:Gsr}. The canonical surjection \(\pi: \bfG(K \otimes_\Q \R) \to \bfG_{s,r}(\R)\) has compact kernel \(\prod_{i=1}^{d-s} \SU_n(\R)\); this implies that the image \(\pi(\Gamma)\) of \(\Gamma\) under \(\pi\) is an irreducible, cocompact (provided \(s < d\)), lattice in \(\bfG_{s,r}(\R)\).

We set \(\gamma\) to be the diagonal \(h\)-unitary \(n \times n\)-matrix \(\diag(\alpha, \alpha^{-1}, 1, \dots, 1)\) in \(\Gamma = \bfG(\cO_K)\). Viewed as an element of \(\bfG_{s,r}(\R)\), the non-trivial block of \(\pi(\gamma)\) is
\[
\left(
\begin{smallpmatrix} \alpha_1 &  \\  & \alpha_1^{-1} \end{smallpmatrix}, \dots,
\begin{smallpmatrix} \alpha_r &  \\  & \alpha_r^{-1} \end{smallpmatrix},
\begin{smallpmatrix} \alpha_{r+1} &  \\  & \alpha_{r+1}^{-1} \end{smallpmatrix}, \dots,
\begin{smallpmatrix} \alpha_{(s+r)/2} &  \\  & \alpha_{(s+r)/2}^{-1} \end{smallpmatrix}
\right),
\]
where we labelled the roots \(\alpha_i\) of \(P\) in \(\C\) as in \S \ref{reduction-subsec}.

\subsection{The sequence \((\gamma'_m)_{m}\) and Margulis' conjecture} \label{sequence-subsec}
Assume \eqref{L} holds, set \(t = (s-r)/2\) and pick a sequence \((P_m)_{m \in \N^*}\) of polynomials in \(\scrP_{s,r}\) with at least one root of absolute value 1, satisfying
\[
1 < \Mahler(P_m) < \exp(\eta_{m,t})\quad
\text{with}\quad\eta_{m,t}=\frac{1}{2m^{t+1}}.
\]
For each \(m\in\N^*\), let us denote
\[
\alpha_m,\quad
L_m,\quad
K_m,\quad
h_m,\quad
\bfG_m,\quad
\Gamma_m,\quad
\gamma_m,\quad\text{and}\quad
\pi_m,
\]
all the objects stemming from the construction in \S\S \ref{numberfields-subsec}--\ref{lattice-subsec} applied to the polynomial \(P_m\).
As above, the non-trivial block of \(\pi_m(\gamma_m) \in \bfG_{s,r}(\R)\) is
\[
\left(
\begin{smallpmatrix} \alpha_{m,1} &  \\  & \alpha_{m,1}^{-1} \end{smallpmatrix}, \dots,
\begin{smallpmatrix} \alpha_{m,r} &  \\  & \alpha_{m,r}^{-1} \end{smallpmatrix},
\begin{smallpmatrix} \alpha_{m,r+1} &  \\  & \alpha_{m,r+1}^{-1} \end{smallpmatrix}, \dots,
\begin{smallpmatrix} \alpha_{m,r+t} &  \\  & \alpha_{m,r+t}^{-1} \end{smallpmatrix}
\right),
\]
where \(\alpha_{m,i}\) are roots of \(P_m\) in \(\C\) labelled according to \S \ref{reduction-subsec}.

Let \(U_m\) denote the inversion-invariant neighborhood of \(1\) in \(\C\) given by
\[
U_m = \left\{z \in \C \,\middle|\, -\frac{1}{m} \leq \log|z| \leq \frac{1}{m} \textrm{ and } -\frac{2\pi}{m} \leq \arg(z) \leq \frac{2\pi}{m}\right\}.
\]
By construction, we have \(1 < |\alpha_{m,i}| \leq \Mahler(P_m) < \exp(\eta_{m,t})\) for \(1 \leq i \leq r+t\). Unfortunately, it may very well be that \(\alpha_{m,i} \notin U_m\). To remedy this, we will use Dirichlet’s simultaneous approximation theorem (see e.g.\ \cite[II, \S 1, Theorem 1A]{Schmidt1980}). Let us recall it under the form we will use.
\begin{lemma*}[Dirichlet]\label{Dirichlet-simultaneous}
Let \(t, m \in \N^*\). For any \(x_1, \dots, x_t \in \R/\Z\), there exists an integer \(0 < c \leq m^t\) for which \(c x_i\) lies in \([-1/m, 1/m] + \Z\) for all \(1 \leq i \leq t\).
\end{lemma*}

\noindent
Applying the approximation theorem to the arguments of \(\alpha_{m,r+1}^2, \dots, \alpha_{m,r+t}^2\) (scaled by \(1/2\pi\)) yields an integer \(c_m\) with \(0 < c_m \leq m^t\), such that
\[\alpha_{m,i}^{2c_m} \in U_{m} \quad \textrm{for} \quad 1 \leq i \leq r+t.\]
Indeed, for \(r+1 \leq i \leq r+t\), by construction of \(c_m\), we have \(-2\pi/m \leq \arg(\alpha_{m,i}^{2c_m}) \leq 2\pi/m\); when \(1 \leq i \leq r\), it is clear that \(\arg(\alpha_{m,i}^2) = 0\). In both cases,
\[
1 < |\alpha_{m,i}^{2c_m}| \leq |\alpha_{m,i}|^{2m^t} < \exp(2m^t\cdot \eta_{m,t})=\exp(1/m).
\]

Finally, let \(\gamma'_m\) denote the element \(\pi_m(\gamma_m^{2c_m}) \in \pi_m(\Gamma_m) \leq \bfG_{s,r}(\R)\). Because \(\{U_m \mid m \in \N^*\}\) forms a basis of neighborhoods of \(1\) in \(\C\), the discussion above (together with the continuity of the regular representation \(L_{\sigma_i}^\times \to \GL_2(K_{\sigma_i})\)) shows that \(\gamma'_m \to 1 \in \bfG_{s,r}(\R)\) as \(m \to \infty\).
In addition, \(\pi_m(\Gamma_m)\) is an irreducible, cocompact lattice in \(\bfG_{s,r}(\R)\) by \S\ref{lattice-subsec}, and, as \(|\alpha_{m,1}| > 1\), \(\gamma'_m\) has infinite order.
This contradicts Conjecture \ref{Margulis-conjecture} for the group \(\bfG_{s,r}\) and concludes the proof of theorem \ref{main-theorem} for \(\bfF = \SL_n\).

\section{The argument for other simple groups}
\label{othergroups-sec}

In this section, we indicate the modifications to run the argument for an isotropic \(\R\)-group \(\bfF\) which is either classical, or adjoint of exceptional type \(\sfG_2\), \(\sfF_4\), \(\sfE_6\), \(\sfE_7\) or \(\sfE_8\). As any absolutely (almost) simple isotropic \(\R\)-group is isogenous to one of these (see e.g.\ \cite{Tits1966}), this is sufficient to conclude the proof of Theorem \ref{main-theorem}. 

For brevity, we omit the parts of the argument which are analogous (if not identical) to their counterparts in \S\S \ref{reduction-subsec}--\ref{sequence-subsec}.
The leitmotiv is the construction of a \(K\)-form \(\bfG\) of \(\bfF\) with the appropriate archimedean shape, in which the norm torus \(\ker \rN_{L/K} \leq \Res_{L/K}(\GL_1)\) embeds over \(\cO_K\).
\medskip

We mostly keep the notation of \S\S \ref{reduction-subsec}--\ref{sequence-subsec}. Given a root \(\alpha\) of some \(P \in \scrP_{s,r}\), the fields \(L\) and \(K\) are constructed identically as in \S \ref{numberfields-subsec}. Let \(L'\) denote the quadratic étale \(K\)-algebra \(K[X] / (X^2+1)\), and let \(\tau'\) denote the non-trivial automorphism of \(L'\) fixing \(K\).
Whenever \(-1\) is not a square in \(K\) (e.g.\ when \(K\) is a real field, which is the case if we pick \(P\) according to \eqref{L}), \(L'\) is a totally complex number field; otherwise, \(L'\) is the product of two copies of \(K\). We will denote by \(\sigma'_i\) one of the two \(\tau'\)-conjugate extensions to \(L'\) of the embedding \(\sigma_i: K \to \C\) (\(1 \leq i \leq d\)).

In addition, we will make use of three quaternion algebras, defined over \(K\) by the following symbols:
\[
D^+ = \< \pantisymmetrized{\alpha}^{2}, \pantisymmetrized{\alpha}^{2} \>_K,\quad
D^- = \< -\pantisymmetrized{\alpha}^{2}, -\pantisymmetrized{\alpha}^{2} \>_K,\quad
D'\: = \< -1, -1 \>_K,
\]
and whose conjugation involutions we denote by \(\tau^+\), \(\tau^-\), and \(\tau'\), respectively. (This last abuse of notation is excused by the fact the restriction of the conjugation involution of \(D'\) to any image of \(L'\) corresponds to \(\tau'\).)
Note that \(\pantisymmetrized{\alpha}^{2}\) does indeed belong to \(K\), and that \(\sigma_i(\pantisymmetrized{\alpha}^{2})\) is positive for \(1 \leq i \leq r\), and negative for \(s+1 \leq i \leq d\).
Thus, the quaternion algebra \(D^+\) (resp.\ \(D^-\), \(D'\)) splits over \(K_{\sigma_i}\) if and only if \(1 \leq i \leq s\) (resp.\ \(r+1 \leq i \leq d\), \(r+1 \leq i \leq s\)). Moreover, as \(L = K(\antisymmetrized{\alpha})\), \(L\) embeds in \(D^+\), or equivalently \(D^+\) splits over \(L\). Similarly, \(L'\) embeds in \(D'\).

Once constructed, we systematically endow the \(K\)-group \(\bfG\) with the \(\cO_K\)-structure obtained by writing the canonical equations defining \(\bfG\) over \(\cO_K\) using the basis \(\{1, \alpha\}\) of \(L\) over \(K\) (and its extension to an \(\cO_K\)-structure on \(D^+\), or the canonical bases of \(L'\), \(D^+\), \(D^-\), \(D'\), depending on the need). Unless indicated otherwise, we then set \(\Gamma = \bfG(\cO_K)\) to be the group of solutions in \(\cO_K\) of these equations. The reader will easily verify that this integral structure has the claimed properties.

For any field extension \(M\) of \(K\) and any \(K\)-algebra \(A\), we set \(A_M = A\otimes_K M\). We also set
\[
\bfG_{s,r} = \prod_{i=1}^r \bfF \times \prod_{i=1}^{(s-r)/2} \Res_{\C/\R}(\bfF) \in \scrT_\bfF^{(s)},
\]
and as in \S\ref{lattice-subsec}, \(\pi\) will be the projection \(\bfG \to \bfG_{s,r}\) with compact kernel to be constructed below. 
\medskip

We now give the modifications by the \(\R\)-type of \(\bfF\), following the descriptions provided by Tits \cite[Table \textrm{II}]{Tits1966} (except perhaps for groups of type \(\sfE\)). The reader may also find \cite{KnusMerkurjevRostTignol1998} a useful reference.
\medskip

\subsection{Inner form of type \(\sfA\).}
\(\bfF\) is of the form \(\SL_{n, D_0}\) for some division algebra \(D_0\) over \(\R\) and \(n \in \N^*\). As the split case has already been covered, we may as well assume that \(D_0 \neq \R\), hence \(D_0\) is simply Hamilton's quaternion algebra. Note that in fact \(n \geq 2\) since \(\bfF\) is isotropic.

Consider the involution \(\tau^- \otimes \tau\) on the quaternion \(L\)-algebra \(D^-_L = D^- \otimes_K L\). As it coincides with \(\tau\) when restricted to the center \(L\) of \(D^-_L\), this involution is of the second kind.
Let \(h: (D^-_L)^n \times (D^-_L)^n \to D^-_L\) be the \(\tau^- \otimes \tau\)-hermitian form given by
\[
h(x,y) = x_1 (\tau^- \otimes \tau)(y_1) + \dots + x_n (\tau^- \otimes \tau)(y_n) \qquad x,\,y\in (D^-_L)^n,
\]
and let \(\bfG\) be the special unitary group \(\SU_{h, D^-_L}\) associated to \(h\).

The local structure of \(G\) is as follows. Of course, \(\bfG\) splits at the complex places \(\sigma_{r+1}, \dots, \sigma_{s}\) of \(K\).
For \(s+1 \leq i \leq d\), the quaternion algebra \(D^-_L\) splits over \(L_{\sigma_i} \cong \C\). Nonetheless, \(\bfG\) remains an outer form (of \(\SL_{2n}\)) over \(K_{\sigma_i}\), since \(\tau^- \otimes \tau\) induces the non-trivial automorphism of the field extension \(L_{\sigma_i} / K_{\sigma_i} \cong \C / \R\). In view of the coefficients of \(h\), we deduce that \(\bfG\) becomes isomorphic to \(\SU_{2n}\) (the usual anisotropic special unitary \(\R\)-group) over \(K_{\sigma_i} \cong \R\) for \(s+1 \leq i \leq d\).
Lastly, for \(1 \leq i \leq r\), we have that \(L\) is contained in \(K_{\sigma_i}\). This means that
\[
D^-_L \otimes_K K_{\sigma_i} \cong D^- \otimes_K (L \otimes_K K_{\sigma_i}) \cong D^- \otimes_K (K_{\sigma_i} \oplus K_{\sigma_i}) \cong D^-_{K_{\sigma_i}} \oplus D^-_{K_{\sigma_i}},
\]
with \(\tau^- \otimes \tau\) inducing a flip of the two summands. Thus \(\bfG \cong \SL_{n, D^-_{K_{\sigma_i}}}\) over \(K_{\sigma_i}\), that is \(\bfG\) becomes an inner form (of \(\SL_{2n}\)) over \(K_{\sigma_i}\).
We had already observed that \(D^-_{K_{\sigma_i}}\) does not split for \(1 \leq i \leq r\), hence \(D^-_{K_{\sigma_i}} \cong D_0\) after identifying \(K_{\sigma_i}\) with \(\R\); in turn, \(\bfG \cong \SL_{n,D_0}\) over \(K_{\sigma_i}\). Altogether, we obtain
\[
\bfG(K \otimes_\Q \R) \cong \prod_{i=1}^{r} \SL_{n,D_0}(\R) \times \prod_{i=1}^{(s-r)/2} \SL_{2n}(\C) \times \prod_{i=1}^{d-s} \SU_{2n}(\R),
\]
so that \(\bfG(K \otimes_\Q \R)\) projects onto \(\bfG_{s,r}(\R)\) with compact kernel.

The element \(\gamma \in \Gamma\) is taken to be \(\diag(\alpha, \alpha^{-1}, 1, \dots, 1)\), seen as an element of \(\bfG(\cO_K) \subset \SL_n(D^-_L)\). One has indeed \(\alpha \cdot \tau(\alpha) = \rN_{L/K}(\alpha) = 1\), so that \(\gamma\) preserves \(h\).

\subsection{Outer form of type \(\sfA\).}
\(\bfF\) is of the form \(\SU_{h_0}\) for some hermitian form \(h_0: \C^n \times \C^n \to \C\) of indefinite signature, say \((p,n-p)\) with \(p \geq n-p\).

The symmetric \(K\)-bilinear form \(b\) on \(L \times L\) associated to the quadratic form \(\rN_{L/K}\) on \(L\) allows us to construct a \(\tau'\)-hermitian form \(h_b\) on \(L \otimes_K L'\) by setting
\[
h_b(l_1 \otimes l'_1, l_2 \otimes l'_2) = b(l_1, l_2) l'_1 \cdot \tau'(l'_2) \qquad l_1,\, l_2 \in L, \; l'_1,\, l'_2 \in L'.
\]
Let \(V\) denote the vector \(L'\)-space \((L \otimes_K L')^{n-p} \oplus L'^{(2p-n)}\), and let \(h: V \times V \to L'\) be the \(\tau'\)-hermitian form given as the orthogonal sum
\[
h=h_{b}^{\oplus (n-p)}\oplus \< 1\>^{\oplus (2p-n)},
\]
where \(\< 1 \>\) denotes the \(\tau'\)-hermitian form \((x,y) \mapsto x \tau'(y)\). Let \(\bfG\) be the special unitary group \(\SU_{h}\) associated to \(h\).

The local structure of \(\bfG\) is as follows. Of course, \(\bfG\) splits at the complex places \(\sigma_{r+1}, \dots, \sigma_s\) of \(K\).
At the other places, one has by construction \(L'_{\sigma'_i}/K_{\sigma_i} \cong \C / \R\) with \(\tau'\) inducing complex conjugation; thus \(\bfG\) remains an outer form (of \(\SL_n\)) there. In addition, the quadratic form \(\rN_{L/K}\) (hence also the hermitian form \(h_b\)) has local signature \((1,1)\) or \((2,0)\) over \(K_{\sigma_i}\) for respectively \(1 \leq i \leq r\) or \(s+1 \leq i \leq d\). Thus \(h\) has local signature respectively \((p,n-p)\) or \((n,0)\) over \(K_{\sigma_i}\).
Altogether, we have that
\[
\bfG(K \otimes_\Q \R) \cong \prod_{i=1}^{r} \SU_{h_0}(\R) \times \prod_{i=1}^{(s-r)/2} \SL_{n}(\C) \times \prod_{i=1}^{d-s} \SU_{n}(\R),
\]
so that \(\bfG(K \otimes_\Q \R)\) projects onto \(\bfG_{s,r}(\R)\) with compact kernel.

Since \(\rN_{L/K}(\alpha)=1\), multiplication by \(\alpha\) induces an isometry of the quadratic \(K\)-space \((L,b)\), and hence a unitary operator \(u_\alpha\) of the hermitian space \((L \otimes_K L', h_b)\). Let then \(\gamma\) be the element of \(\bfG(K)\) which acts as \(u_\alpha\) on the first component of \(V\), and as the identity on the other components. Note that \(\gamma \in \bfG(\cO_K)\) by choice of the \(\cO_K\)-structure.

\subsection{Type \(\sfB\).}
\(\bfF\) is of the form \(\SO_{q_0}\) for some quadratic form \(q_0\) on \(\R^n\) (\(n\) odd) of indefinite signature, say \((p,n-p)\) with \(p \geq n-p\).

Let \(V\) be the \(K\)-module \(L^{n-p} \oplus K^{2p-n}\), and endow \(V\) with the quadratic form \(q\) over \(K\) given as the orthogonal sum
\[
q = \rN_{L/K}^{\oplus (n-p)} \oplus \< 1 \>^{\oplus (2p-n)},
\]
where \(\< 1 \>\) denotes the quadratic form \(x \mapsto x^2\) on \(K\). Let \(\bfG\) be the special orthogonal group \(\SO_q\) associated to \(q\).

\(\bfG\) splits at the complex places \(\sigma_{r+1}, \dots, \sigma_s\) of \(K\). For \(s+1 \leq i \leq d\), the quadratic \(K_{\sigma_i}\)-form \(\rN_{L_{\sigma_i} / K_{\sigma_i}}\) is positive definite (since \(L_{\sigma_i} / K_{\sigma_i} \cong \C / \R\)), whereas for \(1 \leq i \leq r\), it has signature \((1,1)\). Thus \(q\) has local signature respectively \((n,0)\) and \((p,n-p)\). Writing \(\SO_n\) for the usual anisotropic special orthogonal \(\R\)-group, this means
\[
\bfG(K \otimes_\Q \R) \cong \prod_{i=1}^{r} \SO_{q_0}(\R) \times \prod_{i=1}^{(s-r)/2} \SO_{n}(\C) \times \prod_{i=1}^{d-s} \SO_{n}(\R),
\]
so that \(\bfG(K \otimes_\Q \R)\) projects onto \(\bfG_{s,r}(\R)\) with compact kernel.

Since \(\rN_{L/K}(\alpha)=1\), multiplication by \(\alpha\) induces an isometry \(u_\alpha\) of the quadratic \(K\)-space \((L,b)\). We take \(\gamma\) to be the element of \(\bfG(K)\) which acts as \(u_\alpha\) on the first component of \(V\), and as the identity on the other components. Note again that \(\gamma \in \bfG(\cO_K)\) by choice of the \(\cO_K\)-structure.

\subsection{Type \(\sfC\), non-split.}
\(\bfF\) is of the form \(\SU_{h_0, D_0}\), where \(D_0\) is a quaternion \(\R\)-algebra and \(h_0\) is a hermitian form (with respect to quaternion conjugation) on \(D_0^n\), of signature say \((p,n-p)\) with \(p \geq n-p\). Either \(D_0\) is Hamilton's quaternion algebra, or \(D_0\) splits, in which case quaternion conjugation corresponds to adjugation and \(\SU_{h_0, D_0}\) splits to become \(\Sp_{2n}\) regardless of the signature of \(h_0\). We start with the case where \(D_0\) does not split, and thus \(h_0\) is indefinite (because \(\bfF\) is isotropic).

The symmetric \(K\)-bilinear form \(b\) associated to the quadratic form \(\rN_{L/K}\) on \(L\) allows us to construct a \(\tau'\)-hermitian form \(h_b\) on \(L \otimes_K D'\) by setting
\[
h_b(l_1 \otimes d'_1, l_2 \otimes d'_2) = b(l_1, l_2) d'_1 \cdot \tau'(d'_2) \qquad l_1,\, l_2 \in L, \; d'_1,\, d'_2 \in D'.
\]
Let \(V\) be the (left) \(D'\)-module \((L \otimes_K D')^{n-p} \oplus {D'}^{2p-n}\), and endow \(V\) with the \(\tau'\)-hermitian form \(h\) given as the orthogonal sum
\[
h = h_b^{\oplus (n-p)} \oplus \< 1 \>^{\oplus (2p-n)}.
\]
Let \(\bfG\) be the special unitary group \(\SU_{h,D'}\) associated to \(h\).

Of course, \(\bfG\) splits at the complex places \(\sigma_{r+1}, \dots, \sigma_s\). On the other hand, \(D'_{K_{\sigma_i}} \cong D_0\) when \(1 \leq i \leq r\) or \(s+1 \leq i \leq d\), by construction. From the local signatures of \(\rN_{L/K}\), one computes that \(h\) has local signatures \((p,n-p)\) or \((n,0)\) for \(1 \leq i \leq r\) or \(s+1 \leq i \leq d\) respectively. In other words, writing \(\Sp_{n,0}\) for the anisotropic \(\R\)-form of \(\Sp_{2n}\),
\[
\bfG(K \otimes_\Q \R) \cong \prod_{i=1}^{r} \SU_{h_0}(\R) \times \prod_{i=1}^{(s-r)/2} \Sp_{2n}(\C) \times \prod_{i=1}^{d-s} \Sp_{n,0}(\R),
\]
so that \(\bfG(K \otimes_\Q \R)\) projects onto \(\bfG_{s,r}(\R)\) with compact kernel.

Since \(\rN_{L/K}(\alpha)=1\), multiplication by \(\alpha\) induces an isometry of the quadratic \(K\)-space \((L,b)\), and hence a unitary operator \(u_\alpha\) of the hermitian space \((L \otimes_K D', h_b)\). Let then \(\gamma \in \bfG(\cO_K)\) act as \(u_\alpha\) on the first component of \(V\), and as the identity on the other components.

\subsection{Type \(\sfC\), split.}
Now if \(D_0\) is split, we consider instead the \(\tau^+\)-hermitian form \(h = \< 1 \>^{\oplus n}\) on \((D^+)^n\). Since \(D^+\) splits over \(K_{\sigma_i}\) for \(1 \leq i \leq s\), so does the \(K\)-group \(\bfG = \SU_{h,D^+}\). On the other hand, \(D^+\) does not split at the remaining places, and the signature of \(h\) indicates that \(\bfG(K \otimes \R)\) projects appropriately onto \(\bfG_{s,r}(\R)\). It remains to observe that \(\alpha\) can be seen as an element of \(D^+\) through the canonical embedding of \(L\), and viewed as such, the element \(\gamma = \diag(\alpha, 1, \dots, 1) \in \GL_n(D^+)\) actually belongs to \(\bfG(\cO_K)\) because \(\tau^+\) coincides with \(\tau\) on \(L\).

\subsection{Type \(\sfD\).}
\(\bfF\) is either of the form \(\SO_{q_0}\) for some quadratic form \(q_0\) on \(\R^n\) (\(n\) even) of indefinite signature, or of the form \(\SU_{h_0, D_0}\) for \(h_0\) the standard hermitian form on \(D_0^n\) (\(n \geq 2\)) with respect to an involution of orthogonal type, where \(D_0\) is Hamilton's quaternion algebra. The first case is treated exactly like type \(\sfB\); we thus focus on the second case.

Let \(\rho\) denote the involution of orthogonal type on \(D^-\) given by
\[
\rho(x_1 + x_2 \bfi + x_3 \bfj + x_4 \bfk) = x_1 + x_2 \bfi + x_3 \bfj - x_4 \bfk \qquad x_1, x_2, x_3, x_4 \in K
\]
in the canonical basis \(\{1, \bfi, \bfj, \bfk\}\) of \(D^-\) over \(K\). As previously, the symmetric \(K\)-bilinear form \(b\) associated to \(\rN_{L/K}\) allows us to construct a \(\rho\)-hermitian form \(h_b\) on \(L \otimes_K D^-\). 
Let \(V\) be the (left) \(D^-\)-module \((L \otimes_K D^-) \oplus (D^-)^{n-2}\), and endow \(V\) with the \(\rho\)-hermitian form given as the orthogonal sum
\[
h = h_b \oplus \< 1 \>^{\oplus (n-2)}.
\]
Let \(\bfG\) be the special unitary group \(\SU_{h,D^-}\).

By construction, \(D^- \cong D_0\) over \(K_{\sigma_i}\) for \(1 \leq i \leq r\), and \(h\) is easily seen to be equivalent to the standard form \(h_0\) over \(D_0\). For \(s+1 \leq i \leq d\) however, \(D^-\) splits over \(K_{\sigma_i}\), and \(\rho\) becomes the transposition involution on \(\rM_2(K_{\sigma_i})\). Moreover, at these places the form \(b\) is positive definite; so \(h\) is again equivalent to the standard form \(\< 1 \>^{\oplus n}\), showing that over \(K_{\sigma_i}\), \(\bfG \cong \SO_{2n}\) is anisotropic (\(s+1 \leq i \leq d\)). Altogether,
\[
\bfG(K \otimes_\Q \R) \cong \prod_{i=1}^{r} \SU_{h_0,D_0}(\R) \times \prod_{i=1}^{(s-r)/2} \SO_{2n}(\C) \times \prod_{i=1}^{d-s} \SO_{2n}(\R),
\]
so that \(\bfG(K \otimes_\Q \R)\) projects onto \(\bfG_{s,r}(\R)\) with compact kernel.

Again, multiplication by \(\alpha\) induces a unitary operator \(u_\alpha\) on the hermitian space \((L \otimes_K D^-,h_b)\), and we take \(\gamma \in \bfG(\cO_K)\) to act as \(u_\alpha\) on the first component of \(V\), and as the identity on the other components.

\subsection{Type \(\sfG_2\).} \label{G2-subsec}
\(\bfF\) is of the form \(\Aut O_0\), where \(O_0\) is the split octonion \(\R\)-algebra. (Note that regardless of the base field, the only isotropic form of \(\sfG_2\) is the split one.)

Let \(O^+\) be the octonion \(K\)-algebra obtained by applying the Cayley--Dickson construction to \(D^+\) with the parameter \(-1\). Recall that this means \(O^+ = D^+ \oplus D^+ \bfl\), endowed with the multiplication rule
\[
(a+b\bfl)(c+d\bfl) = (ac - \tau^+(d) b) + (da + b\tau^+(c)) \bfl \qquad a, b, c, d \in D^+.
\]
(Recall that a split octonion algebra may be obtained by the Cayley-Dickson construction applied to a split quaternion algebra with any parameter.)

Let \(\bfG\) be the \(K\)-group \(\Aut O^+\).

For \(1 \leq i \leq s\), the algebra \(O^+\), hence also \(\bfG\), splits over \(K_{\sigma_i}\) because \(D^+\) does, whereas \(\bfG\) is anisotropic over \(K_{\sigma_i}\) for \(s+1 \leq i \leq d\), because at these places \(O^+\) is a (nonassociative) division algebra since its norm form (which is the double of the norm form of \(D^+\)) is positive definite. Altogether, writing \(\bfF^\mathrm{anis}\) for the anisotropic \(\R\)-group of type \(\sfG_2\),
\[
\bfG(K \otimes_\Q \R) \cong \prod_{i=1}^{r} \bfF(\R) \times \prod_{i=1}^{(s-r)/2} \bfF(\C) \times \prod_{i=1}^{d-s} \bfF^\mathrm{anis}(\R),
\]
so that \(\bfG(K \otimes_\Q \R)\) projects onto \(\bfG_{s,r}(\R)\) with compact kernel.

One checks that given \(x \in D^+\) with \(\tau^+(x) x = 1\), the assignment \(a + b\bfl \mapsto a + (xb) \bfl\) defines an automorphism \(u_x\) of \(O^+\) (see \cite[\S2.1]{SpringerVeldkamp2000}). As \(L\) embeds in \(D^+\) (with \(\tau^+\) restricting to \(\tau\)), \(\alpha\) can be seen as an element of \(D^+\) of norm 1. The automorphism \(\gamma = u_\alpha\) of \(O^+\) then belongs to \(\bfG(\cO_K)\) (if as before, \(\bfG\) is endowed with the \(\cO_K\)-structure induced by an \(\cO_K\)-structure on \(O^+\) extending the basis \(\{1,\alpha\}\) of \(L \subset D^+\)).

\subsection{Preliminaries for type \(\sfF_4\)} \label{preliminariesF4-subsec}
Recall that given a field \(M\), an octonion algebra \(O\) over \(M\), and parameters \(c=(c_1, c_2, c_3) \in (M^\times)^{3}\), the set
\[
\cA_{M}(O,c) = 
\left \{ \begin{pmatrix}
x_1 & y_3 & c_1^{-1} c_3 \overline{y}_2 \\
c_2^{-1} c_1 \overline{y}_3 & x_2 & y_1 \\
y_2 & c_3^{-1} c_2 \overline{y}_1 & x_3
\end{pmatrix}
\; \middle| \; x_1, x_2, x_3 \in M, y_1, y_2, y_3 \in O \right\}
\]
of \emph{\((c_1, c_2, c_3)\)-hermitian \(3 \times 3\) matrices with entries in \(O\)} (endowed with the Jordan product) forms an exceptional Jordan algebra called an \emph{Albert algebra}. 
When the underlying octonion algebra is split, any choices of \(c_1, c_2, c_3\) yield isomorphic Albert algebras, which are accordingly called \emph{split}. 

Over \(\R\), there are three isomorphism classes of Albert algebras, represented by: the \emph{split} Albert algebra of \((1,1,1)\)-hermitian \(3 \times 3\) matrices with entries in the split octonion \(\R\)-algebra, the algebra of \((1,1,1)\)-hermitian \(3 \times 3\) matrices with entries in Cayley's octonion algebra (which we call the \emph{definite} Albert \(\R\)-algebra), and the algebra of \((1,1,-1)\)-hermitian \(3 \times 3\) matrices with entries in Cayley's octonion algebra (which we call the \emph{indefinite} Albert \(\R\)-algebra). 

\subsection{Type \(\sfF_4\), split} \label{F4split-subsec}
\(\bfF\) is of the form \(\Aut A_0\), where \(A_0\) is the split Albert \(\R\)-algebra. 

Let \(O^+\) be the octonion \(K\)-algebra obtained by applying the Cayley--Dickson construction to \(D^+\) with the parameter \(-1\), and let \(A^+=\cA_{K}\big(O^{+},(1,1,1)\big)\).
Let \(\bfG\) be the \(K\)-group \(\Aut A^+\); it is a simple group of type \(\sfF_4\) (see namely \cite[Ch.~5--7]{SpringerVeldkamp2000} for this fact and others concerning the structure of Albert algebras). 

As mentioned in \S\ref{G2-subsec}, \(O^+\) splits over \(K_{\sigma_i}\) for \(1 \leq i \leq s\), whereas it remains a division algebra over \(K_{\sigma_i}\) for \(s+1 \leq i \leq d\). Thus \(A^+\) hence \(\bfG\) splits over \(K_{\sigma_i}\) for \(1 \leq i \leq s\), while for \(s+1 \leq i \leq d\), \(A^+\) is isomorphic over \(K_{\sigma_i}\) to the definite Albert \(\R\)-algebra, hence \(\bfG\) is isomorphic over \(K_{\sigma_i}\) to the anisotropic simple \(\R\)-group \(\bfF^{\mathrm{anis}}\) of type \(\sfF_4\). Altogether, 
\[
\bfG(K \otimes_\Q \R) \cong \prod_{i=1}^{r} \bfF(\R) \times \prod_{i=1}^{(s-r)/2} \bfF(\C) \times \prod_{i=1}^{d-s} \bfF^\mathrm{anis}(\R),
\]
so that \(\bfG(K \otimes_\Q \R)\) projects onto \(\bfG_{s,r}(\R)\) with compact kernel.

Given \(x \in D^+\) of norm 1, the automorphism \(u_x\) of \(O^+\) from \S\ref{G2-subsec} extends canonically to an automorphism of \(A^+\), still to be denoted \(u_x\). By choice of the \(\cO_K\)-structure on \(L\), and in turn on \(D^+\), \(O^+\) and \(A^+\), the automorphism \(\gamma = u_\alpha\) of \(A^+\) belongs to \(\bfG(\cO_K)\).

\subsection{Type \(\sfF_4\), non-split}
\label{F4nonsplit-subsec}
\(\bfF\) is of the form \(\Aut A_0\), where this time \(A_0\) is the indefinite Albert \(\R\)-algebra.

Let \(O'\) be the octonion \(K\)-algebra obtained by applying the Cayley--Dickson construction to \(D'\) with the parameter \(-1\), and let \(A'=\cA_{K}(O',(1,1,-(\antisymmetrized{\alpha})^{2})\).
Let \(\bfG\) be the \(K\)-group \(\Aut A'\). 

The octonion algebra \(O'\) is isomorphic to Cayley's octonion algebra over \(K_{\sigma_i}\) for all \(1 \leq i \leq d\). The signs of the parameters defining \(A'\) vary under the different embeddings of \(K\) in such a way that for \(1 \leq i \leq s\), \(A'\) is isomorphic over \(K_{\sigma_i}\) to \(A_0\), hence also \(\bfG\) to \(\bfF\). For \(s+1 \leq i \leq d\), \(A'\) is isomorphic over \(K_{\sigma_i}\) to the definite Albert \(\R\)-algebra, hence \(\bfG\) is isomorphic over \(K_{\sigma_i}\) to the anisotropic \(\R\)-form \(\bfF^\mathrm{anis}\), as required. 

Conjugation\footnote{Here, we only use that \(O'\) is a \(K\)-(bi)module, which splits as a direct sum of \(K\) and the \(K\)-submodule of totally imaginary octonions. Conjugation by a matrix with arbitrary entries in \(O'\) makes no sense since \(O'\) is not an associative algebra.} by the matrix
\(
\begin{smallpmatrix}
\myfrac{\psymmetrized{\alpha}}{2} & \0 & \myfrac{\pantisymmetrized{\alpha}^2}{2} \\
\0 & 1 & \0 \\
\myfrac{1}{2} & \0 & \myfrac{\psymmetrized{\alpha}}{2}
\end{smallpmatrix} \in \GL_3(K)
\)
defines a linear map \(u_\alpha \in \GL(A')\). Using the fact that the above matrix preserves the diagonal quadratic form \(\< 1, 1, -\pantisymmetrized{\alpha}^2\>\), it is readily seen that \(u_\alpha\) is in fact an automorphism of the Albert algebra \(A'\). 
What is perhaps less obvious is that \(\gamma = u_\alpha\) belongs to an arithmetic lattice in \(\bfG(K \otimes_\Q \R)\). Let \(\Omega\) be the \(\cO_K\)-submodule of \(\rM_3(O')\) spanned by the canonical basis. Conjugation by the matrix
\(
\begin{smallpmatrix}
1 & \0 & \myfrac{\psymmetrized{\alpha}}{2} \\
\0 & 1 & \0 \\
\0 & \0 & \myfrac{1}{2}
\end{smallpmatrix} \in \GL_3(K)
\)
transforms \(\Omega\) into \(\Omega'\), which has the property that the \(\cO_K\)-submodule \(A' \cap \Omega'\) of \(A'\) is preserved by \(\gamma\). In other words, \(\gamma\) belongs to \(\Gamma = \Stab_{\bfG(K)}(A' \cap \Omega')\), which is a lattice when embedded diagonally in \(\bfG(K \otimes_\Q \R)\). 
This can be seen easily as follows. The quadratic form \(\< 1, -\pantisymmetrized{\alpha}^2 \>\) is nothing but the norm form of \(L/K\) written in the \(K\)-basis \(\{1, \antisymmetrized{\alpha}\}\) of \(L\). 
The matrix
\(
\begin{smallpmatrix}
1 & \myfrac{\psymmetrized{\alpha}}{2} \\
\0 & \myfrac{1}{2}
\end{smallpmatrix}
\)
changes basis back to the standard basis \(\{1, \alpha\}\). The non-trivial block of \(\gamma\) in this new basis then rewrites
\(
\begin{smallpmatrix}
\0 & -1 \\
1 & \symmetrized{\alpha}
\end{smallpmatrix} \in \GL_2(\cO_K),
\) 
which is just the matrix of multiplication by \(\alpha\) on \(L\).\footnote{If there were such things as \emph{hermitian forms over octonions}, one would directly extend the norm form of \(L/K\) to \(O'\) instead of working with the form \(\< 1, -\pantisymmetrized{\alpha}^2\>\). This would simplify the notation and avoid the base-change computations.} This also shows that \(\gamma\) is not torsion, and that \(\gamma^{2c_m} \to 1\) whenever \(\alpha^{2c_m} \to 1\) (as needed in the argument of \S\ref{sequence-subsec}). 

\subsection{Preliminaries for type \(\sfE\)} \label{preliminariesE-subsec}
In order to treat groups of type \(\sfE\), we briefly recall Tits' construction of the exceptional simple Lie algebras. 
Let \(B\) and \(C\) be \emph{composition algebras} over a field \(M\) of characteristic \(\neq 2,3\); this means that \(B, C\) are chosen among \(M\) itself, quadratic étale algebras, quaternion algebras, or octonion algebras over \(M\). Let \(J\) be the Jordan algebra \(\cA_{M}(B,(c_{1},c_{2},c_{3}))\), for some \(c_1, c_2, c_3 \in M\) (defined analogously to \S \ref{preliminariesF4-subsec}). 
The set
\[
\frL(C,J) = \Der(C) \oplus (C^{\circ} \otimes J^{\circ}) \oplus \Der(J),
\]
where \(C^{\circ}\) (resp.\ \(J^{\circ}\)) denotes the kernel of the trace of \(C\) (resp.\ \(J\)), can be endowed with a Lie bracket which turns it into a semisimple Lie \(M\)-algebra whose absolute type is given by the \emph{Freudenthal-Tits magic square}:
\begingroup
\renewcommand{\arraystretch}{1.2}
\[
\begin{array}{|c|c|c|c|c|}
\hline
\hbox{\diagbox{\(\dim_M C\)}{\(\dim_M B\)}} & 1 & 2 & 4 & 8 \\ \hline
1 & \sfA_1 & \sfA_2 & \sfC_3 & \sfF_4 \\ \hline
2 & \sfA_2 & \sfA_2 \times \sfA_2 & \sfA_5 & \sfE_6 \\ \hline
4 & \sfC_3 & \sfA_5 & \sfD_6 & \sfE_7 \\ \hline
8 & \sfF_4 & \sfE_6 & \sfE_7 & \sfE_8 \\ \hline
\end{array}
\]
\endgroup
This construction and the magic square were discovered independently by Freudenthal and Tits. 
Quite remarkably, when extended appropriately to include \(\sfG_2\), the construction gives Lie algebras of all the exceptional types in a unified way. 
For a more detailed description, we refer the reader to the original articles of Freudenthal [\nbcite{Freudenthal1954} -- \nbcite{Freudenthal1963}] \nocite{Freudenthal1954,Freudenthal1954b,Freudenthal1955,Freudenthal1959,Freudenthal1963} and Tits \cite{Tits1966a}, as well as to the lecture notes of Jacobson \cite{Jacobson1971} and the Book of Involutions \cite[Ch.~IX]{KnusMerkurjevRostTignol1998}. 

It turns out that over \(\R\), Tits' construction produces all possible forms of simple real Lie algebras of type \(\sfE\) or \(\sfF\), as \(C\) and \(J\) run through all possible combinations.
We refer to \cite[pp.~119--121]{Jacobson1971} for a summary of which input \((C, J)\) to use for each \(\R\)-form. 
We will use this in the remaining cases (types \(\sfE_6\), \(\sfE_7\), \(\sfE_8\)) to describe the group \(\bfF\) and construct an appropriate form of \(\bfF\) over \(\cO_K\).\footnote{For \(\sfE_6\) (and perhaps for \(\sfE_7\)), one could have used a description more resemblant to that of \(\sfG_2\) and \(\sfF_4\). However, the smallest nontrivial irreducible representation of a group of type \(\sfE_8\) is the adjoint representation, hence one essentially cannot avoid looking at a group of type \(\sfE_8\) as the automorphism groups of its Lie algebra. 
} 
In fact, the case \(\sfF_4\) treated above is also encompassed by this approach (take \(C=\R\) and \(B\) an octonion \(\R\)-algebra); but we decided to treat \(\sfF_4\) beforehand because it will be used for type \(\sfE\), in the following way.

As is easily seen from the definition of the Lie bracket, any automorphism $\varphi$ of \(C\) (resp.\ \(J\)) induces an automorphism $\tilde{\varphi}$ of \(\frL(C,J)\) by acting via conjugation on the Lie subalgebra \(\Der (C)\) (resp.\ \(\Der(J)\)), acting canonically on the first (resp.\ second) tensor component of \(C^{\circ} \otimes J^{\circ}\), and acting trivially on \(\Der(J)\) (resp.\ \(\Der (C)\)). 
This assignment yields morphisms of algebraic groups \(\Aut C \to \Aut \frL(C,J)\), resp.\ \(\Aut J \to \Aut \frL(C,J)\), whose images commute, and which for the last row (resp.\ column) of the magic square correspond to the familiar embeddings \(\sfG_2\) (resp.\ \(\sfF_4\)) \( \hookrightarrow \sfF_4, \sfE_6, \sfE_7, \sfE_8\).

\subsection{Type \(\sfE\)}
\(\bfF\) is of the form \(\Aut \frL(C_0,A_0)\), where \(C_0\) is a quadratic étale, a quaternion or an octonion \(\R\)-algebra (depending on the absolute rank of \(\bfF\)), and \(A_0\) is an Albert \(\R\)-algebra. 
There are respectively four, three and two possibilities for the isotropic \(\R\)-group \(\bfF\). We list them below (labelled by the signature of their Killing form), and fix inputs \((C_0,A_0)\) when there is more than one choice (obviously, a similar construction works with any of the choices). 
We begin with \(\sfE_6\):
\begin{enumerate}[label=(\(\sfE_6\).\roman*),equation]
\item[(\(\sfE_6^{6}\))] \(C_0 \cong \R \oplus \R\) and \(A_0\) is split (as in \S \ref{F4split-subsec}). In this case, \(\bfF\) is split. 
\item[(\(\sfE_6^{-26}\))] \(C_0 \cong \R \oplus \R\) and \(A_0\) is (up to isomorphism) any of the two non-split Albert \(\R\)-algebras; we pick \(A_0\) to be indefinite (as in \S\ref{F4nonsplit-subsec}). In this case, \(\bfF\) is the non-split inner $\R$-form of type \(\sfE_6\). 
\end{enumerate}
The remaining \(\R\)-forms of \(\sfE_6\) are outer forms. The two isotropic ones are given by: 
\begin{enumerate}[label=(\(\sfE_6\).\roman*),leftmargin=*,equation,resume]
\item[(\(\sfE_6^{2}\))] \(C_0 \cong \C\) and \(A_0\) split. 
\item[(\(\sfE_6^{-14}\))] \(C_0 \cong \C\) and \(A_0\) is indefinite (as in \S\ref{F4nonsplit-subsec}). 
\end{enumerate}
Now for \(\sfE_7\):
\begin{enumerate}[label=(\(\sfE_7\).\roman*),leftmargin=*,equation]
\item[(\(\sfE_7^{7}\))] \(C_0\) and \(A_0\) are split. In this case, \(\bfF\) is split. 
\item[(\(\sfE_7^{-25}\))] \(C_0\) is split and \(A_0\) is any of the two non-split Albert \(\R\)-algebras; we pick \(A_0\) indefinite (as in \S\ref{F4nonsplit-subsec}).
\item[(\(\sfE_7^{-5}\))] \(C_0\) is Hamilton's quaternion algebra and \(A_0\) is split or indefinite; we pick \(A_0\) split. 
\end{enumerate}
And lastly for \(\sfE_8\):
\begin{enumerate}[label=(\(\sfE_8\).\roman*),leftmargin=*,equation]
\item[(\(\sfE_8^{8}\))] \(C_0\) and \(A_0\) are split, or alternatively \(C_0\) is Cayley's octonion algebra and \(A_0\) is indefinite; we pick \(C_0\) and \(A_0\) both split. In this case, \(\bfF\) is split. 
\item[(\(\sfE_8^{-24}\))] \(C_0\) is Cayley's octonion algebra and \(A_0\) is split, or alternatively \(C_0\) is split and \(A_0\) is any of the two non-split Albert \(\R\)-algebras; we pick the former. 
\end{enumerate}
Note that for each rank, exactly one combination for \((C_0,A_0)\) is disregarded, namely the one where \(C_0\) is not split and \(A_0\) is definite. This combination gives the (unique) anisotropic \(\R\)-form of type \(\sfE\) of the corresponding rank. 

For each one of these possibilities, we define the \(K\)-algebras \(C\) and \(A\) as indicated in the table below. 
The \(K\)-algebras \(O^+\), \(A^{+}\) are those defined in \S \ref{F4split-subsec}, and \(O'\), \(A'\) are defined in \S \ref{F4nonsplit-subsec}.
\begingroup
\renewcommand{\arraystretch}{1.2}
\[
\begin{tabular}{c|ccccccccc}
\(\R\)-type & \(\sfE_6^{6}\) & \(\sfE_6^{-26}\) & \(\sfE_6^{2}\) & \(\sfE_6^{14}\) & \(\sfE_7^{7}\) & \(\sfE_7^{-25}\) & \(\sfE_7^{-5}\) & \(\sfE_8^{8}\) & \(\sfE_8^{-24}\) \\[1pt]
\hline
\(C\) & \(L\) & \(L\) & \(L'\) & \(L'\) & \(D^{+}\) & \(D^{+}\) & \(D'\) & \(O^{+}\) & \(O'\)\\
\(A\) & \(A^+\) & \(A'\) & \(A^{+}\) & \(A'\) & \(A^{+}\) & \(A'\) & \(A^{+}\) & \(A^{+}\) & \(A^{+}\)
\end{tabular}
\]
\endgroup
We then let \(\bfG\) be the \(K\)-group \(\Aut \frL(C,A)\) afforded by Tits' construction. 

In view of the behavior of the chosen algebras \(C\) and \(A\) at the archimedean places of \(K\) (described in \S\ref{F4split-subsec} or \S\ref{F4nonsplit-subsec}), \(\bfG\) is isomorphic to \(\bfF\) over \(K_{\sigma_i}\) for \(i=1, \dots, s\), whereas it remains anisotropic over \(K_{\sigma_i}\) for \(i=s+1, \dots, d\). 
Thus \(\bfG(K \otimes_\Q \R)\) projects onto \(\bfG_{s,r}(\R)\) with compact kernel. 

Let \(u_\alpha\) be the automorphism of \(A\) constructed from \(\alpha\) in \S\ref{F4split-subsec} (if \(A = A^+\)) or in \S\ref{F4nonsplit-subsec} (if \(A = A'\)). The induced automorphism \(\gamma = \tilde{u}_\alpha\) of \(\frL(C,A)\) is then the required element of \(\bfG(\cO_K)\), where we use for \(\frL(C,A)\) the \(\cO_K\)-structure canonically induced by that of \(C\) and \(A\) (given in \S\ref{F4split-subsec} or \S\ref{F4nonsplit-subsec}). Indeed, \(\gamma\) is not torsion since the morphism \(\tilde{\phantom{u}}: \Aut A \to \Aut \frL(C,A)\) used to construct it is injective, and \(\gamma^{2c_m} \to 1\) whenever \(\alpha^{2c_m} \to 1\) because (the extension of) this morphism (to any completion of \(K\)) is continuous.

\bibliography{bibliography.bib}
\bibliographystyle{alpha}

\end{document}

%% file: math-shortcuts.tex
\usepackage{eucal} 

\newcommand{\frL}{\mathfrak{L}}

\newcommand{\rM}{\mathrm{M}}
\newcommand{\rN}{\mathrm{N}}

\newcommand{\bfC}{\mathbf{C}}

\newcommand{\bfF}{\mathbf{F}}
\newcommand{\bfG}{\mathbf{G}}
\newcommand{\bfH}{\mathbf{H}}

\newcommand{\bfi}{\mathbf{i}}
\newcommand{\bfj}{\mathbf{j}}
\newcommand{\bfk}{\mathbf{k}}
\newcommand{\bfl}{\mathbf{l}}

\newcommand{\cA}{\mathcal{A}}

\newcommand{\cH}{\mathcal{H}}

\newcommand{\cM}{\mathcal{M}}

\newcommand{\cO}{\mathcal{O}}

\usepackage{mathrsfs}

\newcommand{\scrP}{\mathscr{P}}

\newcommand{\scrT}{\mathscr{T}}

\newcommand{\C}{\mathbb{C}}
\newcommand{\N}{\mathbb{N}}
\newcommand{\Q}{\mathbb{Q}}
\newcommand{\R}{\mathbb{R}}
\newcommand{\Z}{\mathbb{Z}}

\newcommand{\ud}{\mathrm{d}}
\DeclareMathOperator{\vol}{vol}

\DeclareMathOperator{\Stab}{Stab}
\DeclareMathOperator{\id}{id}
\DeclareMathOperator{\Gal}{Gal}

\DeclareMathOperator{\Isom}{Isom}

\DeclareMathOperator{\SO}{SO}
\DeclareMathOperator{\GL}{GL}
\DeclareMathOperator{\SL}{SL}

\DeclareMathOperator{\SU}{SU}

\DeclareMathOperator{\Sp}{Sp}

\DeclareMathOperator{\tr}{tr}
\DeclareMathOperator{\rank}{rank}
\DeclareMathOperator{\diag}{diag}

\DeclareMathOperator{\Lie}{Lie}

\DeclareMathOperator{\Aut}{\mathbf{Aut}}
\DeclareMathOperator{\Ad}{Ad}

\DeclareMathOperator{\Der}{Der}

\usepackage{newtxtext}
\usepackage{newtxmath}

\usepackage{etoolbox}
\makeatletter

\ifdef{\makecolors}{
\definecolor{MyDarkGreen}{rgb}{0,0.5,0.1}
\definecolor{MyDarkBlue}{rgb}{0,0.1,0.75}
\definecolor{MyDarkRed}{rgb}{.8,0.1,0}
\definecolor{MyDarkPurple}{rgb}{.33,.1,.55}
}
{\definecolor{MyDarkGreen}{rgb}{0,0,0}
\definecolor{MyDarkBlue}{rgb}{0,0,0}
\definecolor{MyDarkRed}{rgb}{0,0,0}
\definecolor{MyDarkPurple}{rgb}{0,0,0}
}

\makeatother

\newtheoremstyle{mythmstyle}
  {10pt} 
  {10pt} 
  {} 
  {\normalparindent} 
  {\bfseries} 
  {.} 
  {.5em} 
  {} 

\newtheoremstyle{mythmstyle-red}
  {10pt} 
  {10pt} 
  {} 
  {\normalparindent} 
  {\bfseries\color{MyDarkRed}} 
  {.} 
  {.5em} 
  {} 

\newtheoremstyle{mythmstyle-blue}
  {10pt} 
  {10pt} 
  {} 
  {\normalparindent} 
  {\bfseries\color{MyDarkBlue}} 
  {.} 
  {.5em} 
  {} 

\newtheoremstyle{mythmstyle-green}
  {10pt} 
  {10pt} 
  {} 
  {\normalparindent} 
  {\bfseries\color{MyDarkGreen}} 
  {.} 
  {.5em} 
  {} 
  
\newtheoremstyle{mythmstyle-purple}
  {10pt} 
  {10pt} 
  {} 
  {\normalparindent} 
  {\bfseries\color{MyDarkPurple}} 
  {.} 
  {.5em} 
  {} 

\theoremstyle{mythmstyle-red}

\newtheorem*{theorem*}{Theorem}
\newtheorem*{proposition*}{Proposition}
\newtheorem*{corollary*}{Corollary}

\theoremstyle{mythmstyle-purple}

\newtheorem*{lemma*}{Lemma}
\newtheorem*{remark*}{Remark}
\newtheorem*{remarks*}{Remarks}

\theoremstyle{mythmstyle-blue}
\newtheorem{conjecture}{Conjecture}

\newtheorem*{conjecture*}{Conjecture}

\theoremstyle{mythmstyle-green}

\newtheorem*{definition*}{Definition}

\newtheoremstyle{named}{10pt}{10pt}{}{\normalparindent}{\bfseries\color{MyDarkBlue}}{.}{.5em}{\thmnote{#3}}
\theoremstyle{named}
\newtheorem*{statement}{Theorem}

\usepackage[bookmarksdepth=3,colorlinks=true, cite color=red, link color=red]{hyperref}
\usepackage{geometry}
\usepackage{enumitem}

\makeatletter
\renewcommand\subsubsection{\@startsection{subsubsection}{2}%
  \normalparindent{.5\linespacing\@plus.7\linespacing}{-.5em}%
  {\normalfont\bfseries}}
\makeatother

\DeclareSymbolFont{largesymbolsCM}{OMX}{cmex}{m}{n}

\let\prod\relax
\DeclareMathSymbol{\prod}{\mathop}{largesymbolsCM}{"51}

\DeclareSymbolFont{largesymbolsCM}{OMX}{cmex}{m}{n}

\let\coprod\relax
\DeclareMathSymbol{\coprod}{\mathop}{largesymbolsCM}{"60}

\DeclareSymbolFont{largesymbolsCM}{OMX}{cmex}{m}{n}

\let\sum\relax
\DeclareMathSymbol{\sum}{\mathop}{largesymbolsCM}{"50}

\usepackage{bm}
\DeclareMathAlphabet{\mathsfit}{T1}{\sfdefault}{\mddefault}{\sldefault}
\SetMathAlphabet{\mathsfit}{bold}{T1}{\sfdefault}{\bfdefault}{\sldefault}

\newcommand{\sfA}{\mathsf{A}}
\newcommand{\sfB}{\mathsf{B}}
\newcommand{\sfC}{\mathsf{C}}
\newcommand{\sfD}{\mathsf{D}}
\newcommand{\sfE}{\mathsf{E}}
\newcommand{\sfF}{\mathsf{F}}
\newcommand{\sfG}{\mathsf{G}}

\mathcode`l="8000
\begingroup
\makeatletter
\lccode`\~=`\l
\DeclareMathSymbol{\lsb@l}{\mathalpha}{letters}{`l}
\lowercase{\gdef~{\ifnum\the\mathgroup=\m@ne \ell \else \lsb@l \fi}}
\endgroup

\renewcommand{\bfi}{\boldsymbol{i}}
\renewcommand{\bfj}{\boldsymbol{j}}
\renewcommand{\bfk}{\boldsymbol{k}}
\renewcommand{\bfl}{\boldsymbol{l}}

\SetEnumitemKey{equation}{leftmargin=\mathindent,%
	align=left,%
	labelwidth=\mathindent,
	labelsep=0pt}

%% file: main.bbl
\begin{thebibliography}{KMRT98}

\bibitem[BDM07]{BorweinDobrowolskiMossinghoff2007}
Peter Borwein, Edward Dobrowolski, and Michael~J. Mossinghoff.
\newblock Lehmer's problem for polynomials with odd coefficients.
\newblock {\em Ann. of Math. (2)}, 166(2):347--366, 2007.

\bibitem[Ben08]{Benoist2008}
Yves Benoist.
\newblock R\'{e}seaux des groupes de lie.
\newblock Notes de Cours, 2008.

\bibitem[BHC62]{BorelHarish-Chandra1962}
Armand Borel and Harish-Chandra.
\newblock Arithmetic subgroups of algebraic groups.
\newblock {\em Ann. of Math. (2)}, 75:485--535, 1962.

\bibitem[{Bre}07]{Breuillard2007}
Emmanuel {Breuillard}.
\newblock {On uniform exponential growth for solvable groups}.
\newblock {\em {Pure Appl. Math. Q.}}, 3(4):949--667, 2007.

\bibitem[BV20]{BreuillardVarju2020}
Emmanuel Breuillard and P\'{e}ter~P. Varj\'{u}.
\newblock Entropy of {B}ernoulli convolutions and uniform exponential growth
  for linear groups.
\newblock {\em J. Anal. Math.}, 140(2):443--481, 2020.

\bibitem[Dob79]{Dobrowolski1979}
Edward Dobrowolski.
\newblock On a question of {L}ehmer and the number of irreducible factors of a
  polynomial.
\newblock {\em Acta Arith.}, 34(4):391--401, 1979.

\bibitem[ERT19]{EmeryRatcliffeTschantz2019}
Vincent Emery, John~G. Ratcliffe, and Steven~T. Tschantz.
\newblock Salem numbers and arithmetic hyperbolic groups.
\newblock {\em Trans. Amer. Math. Soc.}, 372(1):329--355, 2019.

\bibitem[Fra21]{Fraczyk2021}
Miko\l{}aj Fraczyk.
\newblock Strong limit multiplicity for arithmetic hyperbolic surfaces and
  3-manifolds.
\newblock {\em Invent. Math.}, 224(3):917--985, 2021.

\bibitem[{Fre}54a]{Freudenthal1954}
Hans {Freudenthal}.
\newblock {Beziehungen der \({E_7}\) und \({E_8}\) zur Oktavenebene. I}.
\newblock {\em {Nederl. Akad. Wet., Proc., Ser. A}}, 57:218--230, 1954.

\bibitem[{Fre}54b]{Freudenthal1954b}
Hans {Freudenthal}.
\newblock {Beziehungen der \(E_7\) und \(E_8\) zur Oktavenebene. II}.
\newblock {\em {Nederl. Akad. Wet., Proc., Ser. A}}, 57:363--368, 1954.

\bibitem[{Fre}55]{Freudenthal1955}
Hans {Freudenthal}.
\newblock {Beziehungen der \(E_7\) und \(E_8\) zur Oktavenebene. III, IV}.
\newblock {\em {Nederl. Akad. Wet., Proc., Ser. A}}, 58:151--157, 277--285,
  1955.

\bibitem[{Fre}59]{Freudenthal1959}
Hans {Freudenthal}.
\newblock {Beziehungen der \({\mathfrak E}_ 7\) und \({\mathfrak E}_ 8\) zur
  Oktavenebene. V-IX}.
\newblock {\em {Nederl. Akad. Wet., Proc., Ser. A}}, 62:165--179, 180--191,
  192--201, 447--465, 466--474, 1959.

\bibitem[{Fre}63]{Freudenthal1963}
Hans {Freudenthal}.
\newblock {Beziehungen der \({\mathfrak E}_ 7\) und \({\mathfrak E}_ 8\) zur
  Oktavenebene. X, XL}.
\newblock {\em {Nederl. Akad. Wet., Proc., Ser. A}}, 66:457--471, 472--487,
  1963.

\bibitem[Gel04]{Gelander2004}
Tsachik Gelander.
\newblock Homotopy type and volume of locally symmetric manifolds.
\newblock {\em Duke Math. J.}, 124(3):459--515, 2004.

\bibitem[GH01]{GhateHironaka2001}
Eknath Ghate and Eriko Hironaka.
\newblock The arithmetic and geometry of {S}alem numbers.
\newblock {\em Bull. Amer. Math. Soc. (N.S.)}, 38(3):293--314, 2001.

\bibitem[{Jac}71]{Jacobson1971}
Nathan {Jacobson}.
\newblock {\em {Exceptional Lie algebras}}.
\newblock Lecture Notes in Pure and Applied Mathematics. {New York: Marcel
  Dekker, Inc.}, 1971.

\bibitem[KMRT98]{KnusMerkurjevRostTignol1998}
Max-Albert {Knus}, Alexander {Merkurjev}, Markus {Rost}, and Jean-Pierre
  {Tignol}.
\newblock {\em {The book of involutions. With a preface by J. Tits.}},
  volume~44.
\newblock Providence, RI: American Mathematical Society, 1998.

\bibitem[{Leh}33]{Lehmer1933}
Derrick~H. {Lehmer}.
\newblock {Factorization of certain cyclotomic functions.}
\newblock {\em {Ann. Math. (2)}}, 34:461--479, 1933.

\bibitem[Mar75]{Margulis1975}
Gregory~A. Margulis.
\newblock Discrete groups of motions of manifolds of nonpositive curvature.
\newblock In {\em Proceedings of the {I}nternational {C}ongress of
  {M}athematicians ({V}ancouver, {B}.{C}., 1974), {V}ol. 2}, pages 21--34.
  Canad. Math. Congress, Montreal, Que., 1975.

\bibitem[Mar91]{Margulis1991}
Gregory~A. Margulis.
\newblock {\em Discrete subgroups of semisimple {L}ie groups}, volume~17 of
  {\em Ergebnisse der Mathematik und ihrer Grenzgebiete (3) [Results in
  Mathematics and Related Areas (3)]}.
\newblock Springer-Verlag, Berlin, 1991.

\bibitem[MR03]{MaclachlanReid2003}
Colin Maclachlan and Alan~W. Reid.
\newblock {\em The arithmetic of hyperbolic 3-manifolds}, volume 219 of {\em
  Graduate Texts in Mathematics}.
\newblock Springer-Verlag, New York, 2003.

\bibitem[NR92]{NeumannReid1992}
Walter~D. Neumann and Alan~W. Reid.
\newblock Arithmetic of hyperbolic manifolds.
\newblock In {\em Topology '90 ({C}olumbus, {OH}, 1990)}, volume~1 of {\em Ohio
  State Univ. Math. Res. Inst. Publ.}, pages 273--310. de Gruyter, Berlin,
  1992.

\bibitem[Sch73]{Schinzel1973}
Andrzej Schinzel.
\newblock On the product of the conjugates outside the unit circle of an
  algebraic number.
\newblock {\em Acta Arith.}, 24:385--399, 1973.

\bibitem[Sch80]{Schmidt1980}
Wolfgang~M. Schmidt.
\newblock {\em Diophantine approximation}, volume 785 of {\em Lecture Notes in
  Mathematics}.
\newblock Springer, Berlin, 1980.

\bibitem[Smy71]{Smyth1971}
Chris Smyth.
\newblock On the product of the conjugates outside the unit circle of an
  algebraic integer.
\newblock {\em Bull. London Math. Soc.}, 3:169--175, 1971.

\bibitem[Smy08]{Smyth2008}
Chris Smyth.
\newblock The {M}ahler measure of algebraic numbers: a survey.
\newblock In {\em Number theory and polynomials}, volume 352 of {\em London
  Math. Soc. Lecture Note Ser.}, pages 322--349. Cambridge Univ. Press,
  Cambridge, 2008.

\bibitem[Smy15]{Smyth2015}
Chris Smyth.
\newblock Seventy years of {S}alem numbers.
\newblock {\em Bull. Lond. Math. Soc.}, 47(3):379--395, 2015.

\bibitem[Sur92]{Sury1992}
B.~Sury.
\newblock Arithmetic groups and {S}alem numbers.
\newblock {\em Manuscripta Math.}, 75(1):97--102, 1992.

\bibitem[SV00]{SpringerVeldkamp2000}
Tonny~A. {Springer} and Ferdinand~D. {Veldkamp}.
\newblock {\em {Octonions, Jordan algebras and exceptional groups. Revised
  English version of the original German notes.}}
\newblock Berlin: Springer, revised english version of the original german
  notes edition, 2000.

\bibitem[{Tit}66a]{Tits1966a}
Jacques {Tits}.
\newblock {Alg\`ebres alternatives, alg\`ebres de Jordan et alg\`ebres de Lie
  exceptionnelles. I: Construction}.
\newblock {\em {Nederl. Akad. Wet., Proc., Ser. A}}, 69:223--237, 1966.

\bibitem[Tit66b]{Tits1966}
Jacques Tits.
\newblock Classification of algebraic semisimple groups.
\newblock In {\em Algebraic {G}roups and {D}iscontinuous {S}ubgroups ({P}roc.
  {S}ympos. {P}ure {M}ath., {B}oulder, {C}olo., 1965)}, pages 33--62. Amer.
  Math. Soc., Providence, R.I., 1966.

\bibitem[Vou96]{Voutier1996}
Paul Voutier.
\newblock An effective lower bound for the height of algebraic numbers.
\newblock {\em Acta Arith.}, 74(1):81--95, 1996.

\bibitem[Zim84]{Zimmer1984}
Robert~J. Zimmer.
\newblock {\em Ergodic theory and semisimple groups}, volume~81 of {\em
  Monographs in Mathematics}.
\newblock Birkh\"auser Verlag, Basel, 1984.

\end{thebibliography}
